\documentclass[12pt]{article}
\usepackage{amsmath,amssymb,theorem}
\numberwithin{equation}{section}

 \usepackage{hyperref}

\newtheorem{thm}{Theorem}[section]
\newtheorem{lemma}[thm]{Lemma}
\newtheorem{prop}[thm]{Proposition}
\newtheorem{cor}[thm]{Corollary}
{\theorembodyfont{\rmfamily}

\newtheorem{rmk}[thm]{Remark}}

\renewcommand{\thesubsection}{\arabic{section}.\arabic{subsection}}

\newcommand{\qed}{\hfill \mbox{\raggedright \rule{.07in}{.1in}}}

\newenvironment{proof}{\vspace{1ex}\noindent{\bf
Proof}\hspace{0.5em}}{\hfill\qed\vspace{1ex}}
\newenvironment{pfof}[1]{\vspace{1ex}\noindent{\bf Proof of
#1}\hspace{0.5em}}{\hfill\qed\vspace{1ex}}

\newcommand{\R}{{\mathbb R}}

\newcommand{\Z}{{\mathbb Z}}

\newcommand{\N}{{\mathbb N}}
\newcommand{\T}{{\mathbb T}}

\newcommand{\OO}{{\bf O}}
 \newcommand{\E}{{\bf E}}
 \newcommand{\eps}{{\varepsilon}}
\newcommand{\cA}{{\mathcal A}}
\newcommand{\cK}{{\mathcal K}}
\newcommand{\cL}{{\mathcal L}}
\newcommand{\cM}{{\mathcal M}}
\newcommand{\cO}{{\mathcal O}}
\newcommand{\cP}{{\mathcal P}}

\newcommand{\divv}{\operatorname{div}}

\makeatletter

\title{Quasicrystals in
pattern formation. Part II: Spatially almost-periodic profiles and global existence.}
\author{
Ian Melbourne \thanks{Mathematics Institute,
University of Warwick,
Coventry, CV4 7AL,
UK},
Jens Rademacher \thanks{
Department of Mathematics,
Universit\"at Hamburg,
20146 Hamburg, Germany},
Bob Rink \thanks{
Department of Mathematics,
Vrije Universiteit Amsterdam, The Netherlands},
Sergey Zelik \thanks{Department of Mathematics, ZJNU, Jinhua,  China,\\ and Department of Mathematics, University of Surrey, Guildford, UK\\ and Keldysh institute of applied mathematics, Moscow, Russia\\ and HSE University, Nizgnij Novgorod, Russia}.
}

\date{21 January 2025}

\begin{document}

\maketitle
\begin{flushright}
{\it Dedicated to the fond memory of Claudia Wulff} \hspace*{3em} 
\end{flushright}

 \begin{abstract}
This paper continues our study of quasicrystals initiated in Part~I.  We propose a general mechanism for  constructing quasicrystals, existing globally in time, in spatially-extended systems (partial differential equations with Euclidean symmetry) and demonstrate it on model examples of the Swift-Hohenberg and Brusselator equations. In contrast to Part~I,  our approach here emphasises the theory of almost-periodic functions as well as the global solvability of the corresponding equations in classes of spatially non-decaying functions.
We note that the existence of  such time-evolving quasicrystals with rotational symmetry of all orders, icosahedral symmetry, etc., does not require  technical issues such as Diophantine properties and hard implicit function theorems,  which look unavoidable in the case of steady-state quasicrystals.

This paper can be largely read independently of Part~I. Background material and definitions are repeated for convenience, but some elementary calculations from Part~I are omitted.
 \end{abstract}

\tableofcontents
 \section{Introduction}

The existence of quasicrystals~\cite{SBGC} was first reported in 1984, and has been the source of much interest ever since. Subsequently, there have been several instances of quasicrystalline solutions in fluid experiments. For example,
quasipatterns with eightfold symmetry~\cite{CAL} and
twelvefold symmetry~\cite{EF} were observed in the Faraday wave experiment and quasipatterns with twelvefold symmetry were observed in shaken convection~\cite{VM}. Figures from~\cite{SBGC,EF} are shown in the first part~\cite{MRRZ1} of this paper which is from now on referred to as Part~I.

It is well-known that in systems of partial differential equations (PDEs) with Euclidean $\E(d)$ symmetry, the variation of a parameter generically gives rise to a large variety of spatially-periodic solutions~\cite{BCM,CrawKnob91,CrossHoh,DG,GSS88,GS02,Michel80,Satt79}.
This mechanism is known as \emph{spontaneous symmetry breaking}.
In Part~I, we pointed out that such bifurcations automatically give rise also to a large class of quasicrystals.

The quasicrystals $u:\R^d\to\R$ considered in Part~I have the desirable properties that
\begin{itemize}
\item[(i)]
There is a finitely-generated relatively dense  subgroup $\cL^*\subset\R^d$ such that
$u(x)=\sum_{k\in\cL^*}a_ke^{ik\cdot x}$ with
amplitudes
$a_k\in\R$ satisfying $\sum_k |a_k|<\infty$;
\item[(ii)] The subgroup of $\cL^*$ generated by $\Lambda_{u,\eps}=\{k\in\cL^*:|a_k|>\eps\}$ is not uniformly discrete for some $\eps>0$.
\end{itemize}
Often we are also able to verify
the following condition which strengthens condition~(ii):
\begin{itemize}
\item[(iii)] For any $M,r>0$ there exists $\eps>0$ such that
$\Lambda_{u,\eps}$ is $r$-dense in $\cL^* \cap B(0,M)$.
\end{itemize}
See Part~I for a comparison of these notions with mathematical quasicrystals.

Recall that $\E(d)$ can be regarded as a semidirect product of
$\OO(d)$ (rotations and reflections) and $\R^d$ (translations).
Let $H$ be a finite subgroup of $\OO(d)$ and fix a unit vector
$k_0\in\R^d$. Define $\cL^*_H$ to be the subgroup of $\R^d$  generated by
the vectors $\gamma k_0$, $\gamma\in H$. We require further that $H$ is the largest subgroup of $\OO(d)$ that preserves $\cL^*_H$ so $H$ is a \emph{holohedry} (in particular, $-I\in H$).
We then restrict to functions $u$ that are $H$-invariant, so
the amplitudes $a_k$ satisfy $a_{\gamma k}=a_k$ for all $\gamma\in H$, $k\in\cL^*_H$.

As pointed out in Section~2.3 of Part~I,  the local existence of  solutions of the type 
$$
u(x)=\sum_{k\in\cL^*_H}a_k e^{ik\cdot x},
\qquad  a_{\gamma k}=a_k
\quad\text{for all $\gamma\in H$, $k\in\cL^*_H$},
$$
with $\sum_{k\in\cL^*_H}|a_k|<\infty$
are guaranteed via spontaneous symmetry breaking. In cases where the \emph{holohedry} $H$ violates the \emph{crystallographic restriction} (eg.\ $q$-fold rotations with $q\ge 8$ even or icosahedral symmetry in $\R^3$) these solutions are quasicrystals satisfying conditions~(i) and~(ii) above. 
In this generality, the quasicrystal solutions are guaranteed to exist for a finite amount of time, but may blow up in finite time.

Here, we are mainly interested in the global time existence of such solutions. 
For this, it turns out to be be natural to relax the absolute summability condition in~(i) and to consider classes of almost periodic functions.
We suggest a general approach based on the theory of spatially almost-periodic functions and global well-posedness results for spatially extended systems in classes of spatially non-decaying solutions.

The  paper is organised as follows.
 In Section \ref{sec:AP}, we discuss known facts from the theory of almost-periodic and quasiperiodic functions, which are crucial in what follows.
In Section~\ref{sec:SH}, we illustrate our  approach with results on existence of global branches of quasicrystals in  a relative simple example, the Swift-Hohenberg equation~\cite{SwiftHoh}.  
 Our treatment includes standard facts about the existence of a global attractor, as well as more delicate issues related to spatially almost-periodic solutions.
In Section~\ref{sec:SH2}, we obtain further results that rely more on the specific nature of the Swift-Hohenberg equation.
Finally, in Section \ref{sec:beyond}, we briefly discuss how our approach works for another model example of a reaction-diffusion system with a Turing instability, namely, the Brusselator model.

\section{Almost-periodic functions and quasicrystals}\label{sec:AP}
In this section, we introduce some technical tools and definitions which will be used throughout the paper. We start with spatially almost-periodic and quasiperiodic functions in Subsection~\ref{sec:ap} before specialising to quasicrystals in Subsection~\ref{sec:qc}.

\subsection{Almost-periodic and quasiperiodic functions}\label{sec:ap}

Recall that a set is called \emph{relatively dense} if there is an $R>0$ so that the set intersects every ball of radius $R$.
A continuous function $u:\R^d \to \R$ is called \emph{(spatially) almost-periodic}  if for every $\eps>0$ there is a relatively dense set of vectors $h \in \R^d$ such that 
$$\| u(\cdot + h) - u(\cdot) \|_{\infty} = \sup_{x\in \R^d} |u(x+h)-u(x)| < \eps .
$$
A vector $h$ satisfying this inequality is called an $\eps$-almost-period of $u$.  
The collection of almost-periodic functions on $\R^d$ will be denoted $AP(\R^d)$. Every almost-periodic function is uniformly continuous and
bounded, so $AP(\R^d) \subset C_b(\R^d)$. Sums, products, and  uniform limits of almost-periodic functions are again almost-periodic; hence $AP(\R^d)$ is a Banach-subalgebra of $C_b(\R^d)$.

Given $u\in AP(\R^d)$, the limit $a_u(k)=\lim_{M\to\infty}\frac{1}{(2M)^d}\int_{[-M,M]^d} u(x)e^{-ik\cdot x}\,dx$ exists for all $k\in\R^d$ and the set
$F_u=\{k\in\R^d:a_u(k)\neq0\}$
is at most countable.
A countable subgroup $\cL^*\subset\R^d$ is called a \emph{frequency module} of $u$ if $F_u\subset\cL^*$ and is called \emph{the} frequency module if
$\cL^*$ is generated by $F_u$.

Let $\cL^*\subset\R^d$ be a countable subgroup. We define $AP(\cL^*)$ to be the vector space of almost-periodic functions for which $\cL^*$ is a frequency module. Functions $u\in AP(\cL^*)$ have the Fourier series representation
\begin{equation}\label{usum}
u(x) \sim \sum_{k\in \mathcal{L}^*} a_k e^{i k \cdot x}
\quad\text{where $a_k=a_u(k)$}.
\end{equation}

A classical result of Bohr states that $u\in AP(\R^d)$ if and only if $u$ is a uniform limit of trigonometric polynomials, in which case it is a uniform limit of trigonometric polynomials contained in $AP(\cL^*)$ where $\cL^*$ is the frequency module for $u$.

We define $\ell^1(\cL^*)$ and $\ell^2(\cL^*)$ to consist of functions of the form \eqref{usum} for which the amplitudes $a_k$ are absolutely summable, respectively square summable, with norms
$$
\|u\|_{\ell^1}:=\sum_{k\in\cL^*}|a_k|,\qquad  \|u\|_{\ell^2}^2:=\sum_{k\in\cL^*}|a_k|^2.
$$
It is immediate from Bohr's result that $\ell^1(\cL^*)\subset AP(\cL^*)$ and that $\|u\|_{\infty} \leq \|u\|_{\ell^1}$. 
The Besicovitch space $\ell^2(\cL^*)$ can be defined as the completion of $AP(\cL^*)$ with respect to the $\ell^2$ norm; hence
$$
\ell^1(\cL^*)\subset AP(\cL^*)\subset \ell^2(\cL^*)\cap C_b(\R^d).
$$
Moreover, it follows from Parseval's identity that
$$
\|u\|^2_{\ell^2(\cL^*)}= \lim_{M\to\infty}\frac{1}{(2M)^d}\int_{[-M,M]^d} u(x)^2\,dx \quad\text{for all $u\in AP(\cL^*)$}.
$$
This shows in particular that $\|u\|_{\ell^2} \leq \|u\|_{\infty}$ for $u\in AP(\R^d)$.

We refer to \cite{Bes54,LeZhi} for further details and more general classes of functions generated by trigonometric polynomials.

If $\cL^*$ is finitely generated over $\Z$, then functions $u\in AP(\cL^*)$ are called \emph{quasiperiodic}. 
Assume that $\cL^*$ is generated by the frequency vectors $k_1, \dots, k_p$. This time, we require that these generators are linearly independent over $\Z$, so that the map
$$
m \mapsto m_1 k_1 + \dots + m_p k_p
$$
is a bijection from $\Z^p$ to $\cL^*$. 
 Define 
\begin{equation}\label{Adefinition}
A:\R^d\to\T^p, \qquad Ax = (k_1\cdot x, \dots, k_p\cdot x) \!\!\! \mod (2\pi\Z)^p .
\end{equation}
By linear independence of $k_1, \dots, k_p$ over $\Z$, the image of $A$ is dense in $\T^p$. 

To each quasiperiodic function $u\in AP(\cL^*)$, we associate the \emph{hull function} $U : \T^p \to \R$  defined by
\begin{equation} \label{eq:hull}
U(\phi) := \sum_{m\in  \Z^p} a_{{m_1 k_1 + \dots + m_p k_p}} e^{i m \cdot \phi},
\quad \phi\in\T^p.
\end{equation}
With this definition, we find that
$$
U(Ax)=U(k_1\cdot x, \dots, k_p\cdot x)  = \sum_{m\in  \Z^p} a_{m_1 k_1 + \dots +m_p k_p} e^{i (m_1 k_1 +\dots + m_p k_p) \cdot x}  = \sum_{k\in \cL^*} a_k e^{i k \cdot x} = u(x) .
$$
Since the image of $A$ is dense, $U$ must be a uniform limit of trigonometric polynomials (because $u$ is). Hence, $U\in C(\T^p)\subset L^2(\T^p)$. It is also clear that $\|U\|_{\infty}= \sup_{\phi\in\T^p}|U(\phi)|  = \sup_{x\in \R^d}|u(x)|  = \|u\|_{\infty}$. 
Moreover, for $u, v \in AP(\cL^*)$, the
$L^2$-inner product $\langle U, V\rangle_{L^2} := \frac{1}{(2\pi)^p} \int_{\T^p}  U(\phi) \overline{V}(\phi)\,d\phi$ of their hull functions $U,V$ satisfies
$$
\langle U, V\rangle_{L^2}   = \sum_{m\in \Z^p} a_{m_1 k_1 + \dots +m_p k_p} \overline{b_{m_1 k_1 + \dots +m_p k_p}}  = \sum_{k\in \cL^*} a_k \overline{b_k} = \langle u, v \rangle_{\ell^2}\ ,
$$
i.e., the $\ell^2$-inner product of two quasiperiodic functions equals the $L^2$-inner product of their hull functions. In particular, $\| u\|_{\ell^2} = \|U\|_{L^2}$.

Conversely, if $u: \R^d\to \R$ satisfies $u(x)=U(Ax)$ for some continuous function $U:\T^p\to \R$ and an $A: \R^d \to \T^p$ as given in \eqref{Adefinition} having a dense image, then $u \in AP(\cL^*)$, where $\cL^*$ is the module generated by the rows of $A$, and $U$ is its hull function. As an immediate consequence, we see that $AP(\cL^*)$ is a Banach algebra.

 We finish this section with an alternative characterisation of $AP(\cL^*)$ that we need to prove that well-posedness of a PDE in $C_b(\R^d)$ implies  well-posedness in $AP(\cL^*)$.

\begin{prop}\label{characterisation}
Let $\cL^*$ be finitely generated by $\Z$-linearly independent vectors $k_1, \dots, k_p$ and define $A:\R^d\to\T^p$ as in \eqref{Adefinition}. 
Suppose that $u\in C_b(\R^d)$. 
Then the following are equivalent.
\begin{itemize}
\item[(i)]  $u\in AP(\cL^*)$; and
\item[(ii)] For every sequence $h_j\in \R^d$ satisfying $\lim_{j\to\infty} Ah_j= 0$,  the difference $u(\cdot + h_j) - u(\cdot)$ converges  to zero uniformly. 
\end{itemize}
\end{prop}

\begin{proof}
 (i) $\implies$ (ii). 
Let $h_j$ be a sequence with $Ah_j\to 0$ and let $\eps>0$. As the hull function $U$ of $u$ is uniformly continuous, there exists a $\delta>0$ so that $|\phi - \psi|<\delta$ implies $|U(\phi)-U(\psi)| < \eps$. Let $N$ be so large that $|Ah_j| < \delta$ for all $j>N$. Then $|u(x+h_j) - u(x)| = |U(Ax+Ah_j) - U(Ax)| < \eps$ for $j>N$, so $u(\cdot + h_j)-u(\cdot)$ converges to zero uniformly.

 (ii) $\implies$ (i). 
 We start by showing that a function $U:\T^p\to\R$ can be defined by $U(\phi) := \lim_{j\to \infty}  u(S_j)$ for any sequence $S_j\in \R^d$ with $\lim_{j\to\infty} AS_j = \phi$.
 First of all, note that the sequence $u(S_j)\in \R$ is bounded because $u\in C_b(\R^d)$. Hence it has a converging subsequence, i.e., there exists at least one sequence $S_j$ with $AS_j\to \phi$ and for which $\lim_{j\to\infty} u(S_j)$ exists. If $s_j$ is another such sequence, then $A(S_j-s_j) \to 0$ and therefore $u(x+ S_j-s_j) - u(x) \to 0$ uniformly in $x$ by assumption. 
Choosing $x=s_j$ yields that $\lim_{j\to\infty} u(S_j) =\lim_{j\to\infty} u(s_j)$.  Thus, $U$ is well-defined.

Choosing $S_j=x$ in the definition of $U$ we see that $u(x)=U(Ax)$. 
It remains to show that $U$ is continuous: this implies that $u\in AP(\cL^*)$ and has $U$ as its hull function. To this end, choose $\phi \in \T^p$, a sequence $\phi_j \to \phi$,  a sequence $S_j\in \R^d$ with $\lim_{j\to\infty} AS_j=  \phi$,  and
 sequences $S^j_n$ with $\lim_{n\to\infty} AS^j_n =  \phi_j$. Let $\eps>0$ be given.
The uniform continuity of $u$ implies that there is a $\delta>0$ such that $\|S_m - S_n^j\| < \delta$ implies $\|u(S_m)-u(S_n^j)\| < \eps/3$. We now choose $j, m, n$ large enough so that $\|S_m-S^j_n\|<\delta$ and $\|U(\phi)-u(S_m)\|<\eps/3$ and $\|U(\phi_j)-u(S_n^j)\|<\eps/3$. Then
$$\| U(\phi) - U(\phi_j)\| \leq \| U(\phi) - u(S_m)\| + \|u(S_m) - u(S^j_n)\| + \|u(S^j_n) - U(\phi_j)\| < \eps \, .$$
This proves that $U(\phi_j)\to U(\phi)$ and hence $U$ is continuous.
\end{proof}

\subsection{Quasicrystals}\label{sec:qc}

We now focus on quasicrystals, which we define as a specific type of quasiperiodic function $u\in AP(\R^d)$. 
Let $H$ be a finite subgroup of $\OO(d)$  and fix a vector $k_0\in\R^d$.
Define
$\cL^*_H$ to be the subgroup of $\R^d$  generated by $\{\gamma k_0:\gamma\in H\}$. 
We require that $\cL^*_H$ is not uniformly discrete.
Also, we assume that $H$ is the maximal subgroup of $\OO(d)$ preserving $\cL^*_H$ so $H$ is the \emph{holohedry} of $\cL^*$; in particular $-I\in H$.

We consider quasiperiodic functions $u:\R^d\to\R$ with Fourier series representation
$u(x)\sim \sum_{k\in\cL^*_H}a_ke^{ik\cdot x}$ where $a_k\in\R$.
As in Part~I, by a slight abuse of notation we define $AP(\cL^*_H)$ to consist of such functions satisfying the additional constraint
$a_{k} = a_{\gamma k}$ for all $\gamma \in H$.
(Since we are interested in real-valued functions and $-I\in H$, no extra generality is gained by allowing complex amplitudes $a_k$.)
Similarly, we define $\ell^1(\cL^*_H)$ and $\ell^2(\cL^*_H)$.

\begin{rmk}
A standard calculation, see Proposition~2.4 of Part~I, shows that
$u(\gamma x)=u(x)$ for all $u\in \ell^1(\cL_H^*)$, $\gamma\in H$,
$x\in\R^d$.
The same argument works for $u\in\ell^2(\cL_H^*)$ and hence for $u\in AP(\cL_H^*)$.
\end{rmk}

We say that a function $u:\R^d\to \R$ is \emph{spatially constant} if $u(x)\equiv c$ for some $c\in\R$.
Since $\cL^*_H$ is not uniformly discrete, any non-spatially constant element of $AP(\cL^*_H)$ is automatically a quasicrystal satisfying conditions~(i) and~(ii)  from the introduction. 
Typical functions in $AP(\cL^*_H)$ satisfy condition~(iii).
In the case of time-dependent quasicrystals, it is desirable that~(iii) is satisfied for almost all time. This property is verified where possible in this paper.

\section{The Swift-Hohenberg equation; general arguments}\label{sec:SH}

In this section, we discuss general features related to global existence of spatially-extended PDEs.
For definiteness, we consider 
a relatively simple example, namely the Swift-Hohenberg equation~\cite{SwiftHoh} given by
\begin{align} \label{eq:SH}
\partial_tu=F(u,\lambda)=-(\Delta+1)^2u+\lambda u-u^3, \qquad  u|_{t=0}=u_0.
\end{align}
Here, $u=u(t,x)=u(t)(x):[0,\infty)\times\R^d\to\R$ is an unknown function, $\lambda\in\R$ is a parameter and $u_0=u_0(x)$ is the initial data.
Equation \eqref{eq:SH} undergoes a bifurcation as $\lambda$ passes through zero, with nonzero critical wavenumber $k_c=1$. That is,
the linearisation $(dF)_{0,0}$ has a zero eigenvalue with kernel consisting of all wavefunctions $e^{ik\cdot x}$ with $|k|=1$.
This is called a \emph{bifurcation of type \emph{I}$_s$} in~\cite{CrossHoh} and is called a
\emph{steady-state bifurcation with nonzero critical wave number} in~\cite{M99}.

In Subsection \ref{sec:globalL8}, we recall known facts related with the global well-posedness, regularity and smoothing property of bounded solutions of the Swift-Hohenberg equation in $\R^d$, $d\le 9$, and discuss the dissipativity of the associated solution semigroup $S_\lambda(t):L^\infty(\R^d)\to L^\infty(\R^d)$ as well as the corresponding global attractors $\cA_\lambda$. 
We assume throughout that $d\le 9$; for $d\ge 10$, it is not clear how to obtain global existence.

Although most of these results are well-known, we prefer to give a brief exposition since it forms the basis of our study of equation~\eqref{eq:SH} in the phase spaces related with almost-periodic functions and quasicrystals. Details of the proof of global existence in $L^\infty(\R^d)$ are given in Appendix~\ref{app:global}. 

In Subsection \ref{sec:globalAP}, we begin our study of phase spaces of spatially almost-periodic functions. Using standard arguments based on the global well-posedness in $L^\infty(\R^d)$, we show that if the initial data is almost-periodic, then the corresponding solution remains almost-periodic with the same frequency module for all times $t\ge0$. 
The $\ell^2$ norm of these solutions is automatically bounded in time (since this is true in $L^\infty(\R^d)$).

In Subsection~\ref{sec:globall1}, we extend the approach of \cite{GigaMahalovYoneda11} to construct a wide class of spatially quasiperiodic solutions with finite $\ell^1$ norm for all $t\ge0$, see Theorem \ref{thm:Hs}.
The $\ell^1$ norm grows at most exponentially as $t\to\infty$ but may be unbounded.

Specialising to frequency modules $\cL^*_H$ that are not uniformly discrete, it is an immediate consequence of these results that
we have global existence of quasicrystal solutions with bounded $\ell^2$ norm in each such $AP(\cL^*_H)$ and that a wide subclass of such solutions have finite $\ell^1$ norm.

\subsection{Well-posedness, dissipativity and smoothing in $L^\infty(\R^d)$}\label{sec:globalL8} 

We work throughout with solutions to equation~\eqref{eq:SH} that are spatially bounded.
It is more standard in the literature to consider uniformly local spaces $L^2_b(\R^d)$ but working in $L^\infty(\R^d)$ enables the application of the arguments 
in Section~\ref{sec:ap}.

Let us assume that $u_0\in L^\infty(\R^d)$. Then a function $u\in L^\infty([0,T];L^\infty(\R^d))\cap C([0,T];L^2_{loc}(\R^d))$ is a solution of equation \eqref{eq:SH} if it satisfies the equation in the sense of distributions.
 
For moderate values of $d$ (specifically $d\le 9$), for any $u_0\in L^\infty(\R^d)$, equation~\eqref{eq:SH} possesses a unique global solution $u(t)\in L^\infty(\R^d)$,  $t>0$, see e.g.\ \cite{MiranvilleZelik} or \cite{efendiev}. 
This solution remains bounded as $t\to\infty$
  and the equation generates a dissipative solution semigroup $S_\lambda(t)$, $t\ge0$, in the phase space $L^\infty(\R^d)$: for any $\lambda_0>0$,
there exist a constant $C\in\R$ and a monotone increasing function $Q$ such that
$$
\|S_\lambda(t)u_0\|_{L^\infty}\le Q(\|u_0\|_{L^\infty})e^{-t}+C,\quad t\ge0,
$$
for all $u_0\in L^\infty(\R^d)$ and
$\lambda\in(-\infty,\lambda_0]$. 
Such results are well-known to experts and we provide
details of the arguments in Appendix~\ref{app:global}. 

Moreover \cite{TakTi}, we have \emph{parabolic smoothing}: solutions become analytic in space and time for $t>0$ and bounded in a strip:
$$
|u(t+i\eta,x+iy)|\le C_\eps,\ \ t\ge\eps>0,\ \ |\eta|+|y|\le \kappa_\eps>0.
$$

Combining dissipativity and parabolic smoothing,
 we conclude that the semigroup $S_\lambda(t)$ possesses a locally compact global attractor $\cA_\lambda$, that is:
\begin{itemize}
\item $\cA_\lambda$ is a bounded closed subset of $L^\infty(\R^d)$, which is compact in $L^\infty_{loc}(\R^d)$; 
\item $\cA_\lambda$ is strictly invariant: $S_\lambda(t)\cA_\lambda=\cA_\lambda$ for all $t\in\R$; 
\item $\cA_\lambda$ attracts the images of bounded sets of $L^\infty(\R^d)$ in the topology of $L^\infty_{loc}(\R^d)$, i.e., for every $B$ bounded in $L^\infty(\R^d)$ and every neighbourhood $\cO(\cA_\lambda)$ of the set $\cA_\lambda$ in the topology of $L^\infty_{loc}(\R^d)$, there exists $T=T(\cO,B)$ such that
$$
S_\lambda(t)B\subset\cO(\cA_\lambda) \quad\text{for all  } t\ge T.
$$
\end{itemize}
It follows from the general theory that the family of attractors $\cA_\lambda$ is upper semicontinuous with respect to $\lambda\in(-\infty,\lambda_0]$, see \cite{MiranvilleZelik,Zel23}.  
Moreover, 
$$
\cA_\lambda=\cK_{\lambda}\big|_{t=0},
$$
where $\cK_\lambda\subset L^\infty(\R\times\R^d)$ is the set of all complete bounded solutions of~\eqref{eq:SH}.
 
Let us now discuss continuity of solutions with respect to initial data. 
The standard $L^\infty$ estimate applied to the linear equation for the difference of two solutions $u(t)$ and $v(t)$ gives
\begin{equation}\label{dif}
\|u(t)-v(t)\|_{\infty}\le Ce^{Kt}\|u_0-v_0\|_{\infty},
\end{equation}
where the constants $C$ and $K$ may depend on the $L^\infty$ norms of $u_0$ and $v_0$, but are independent of $t$. Combining this estimate with the parabolic smoothing property, we end up with
\begin{equation}\label{smo}
\|u(t)-v(t)\|_{W^{s,\infty}}\le C_s t^{-N_s}e^{Kt}\|u_0-v_0\|_{L^\infty},\ t>0,
\end{equation}
for $s>0$, where the constants $N_s$, $C_s$ and $K$ are independent of time, see e.g.\ \cite{efendiev}.
 
Continuity with respect to time is a bit more delicate. Indeed, although any $C_0$-semigroup in $L^\infty$ is generated by a bounded operator, even if we replace $L^\infty(\R^d)$ by $C_b(\R^d)$, we may not have continuity in time, see e.g.\ \cite{MiranvilleZelik}. 
However, if we further restrict the phase space and consider the subspace $\widetilde C_b(\R^d)$ given by the closure of $C^1_b(\R^d)$ in $C_b(\R^d)$, then the linear part of equation \eqref{eq:SH} generates an analytic
semigroup in $\widetilde C_b(\R^d)$, see \cite{heck}. Therefore the trajectory $t\mapsto u(t)$ starting from $u_0\in\widetilde C_b(\R^d)$ is continuous with respect to time at $t=0$ in the $L^\infty$ norm. 
 
\subsection{Global existence and boundedness in $AP(\cL^*)$}\label{sec:globalAP}

Global well-posedness (i.e.\  existence and uniqueness, continuous dependence on initial conditions) in the space $L^\infty(\R^d)$ of bounded functions is known for  numerous $\E(d)$-equivariant PDEs, especially if the linear part of the PDE possesses certain parabolic smoothing properties. 
As explained below, this in turn implies boundedness and well-posedness in $AP(\cL^*)$.
In particular, our results apply to the Swift-Hohenberg equation~\eqref{eq:SH} provided $d\le 9$.
 Arguing in a standard way, we get the following result.

 \begin{prop}\label{prop:globalAP} Let $u_0\in AP(\cL^*)$ for some countable subgroup $\cL^*\subset\R^d$. Then $u(t)\in AP(\cL^*)$ for all $t>0$ and $\sup_{t>0}\|u(t)\|_{\ell^2}<\infty$. Moreover, $D^{k}_tD^{s}_x u(t)\in AP(\cL^*)$ for all $k,s\in\N$ and $t>0$.
 \end{prop}

 \begin{proof} Let $u_0\in AP(\cL^*)$. Then using  estimate \eqref{dif} with $u_0(x)$ and $u_0(x+h)$ for any $h\in\R^d$ and utilising the spatial homogeneity of equation~\eqref{eq:SH}, we see that
 $$
 \|u(t,\cdot)-u(t,\cdot+h)\|_{L^\infty}\le Ce^{Kt}\|u_0(\cdot)-u_0(\cdot+h)\|_{L^\infty},
 $$
 for any $t>0$ fixed and any $h\in\R^d$. Thus, any $\eps$-almost period $h$ of $u_0$ gives us a
 $\delta:=Ce^{Kt}\eps$-almost period of $u(t)$ and therefore $u(t)$ is also almost-periodic for every fixed $t\ge0$.  
Moreover, $\|u(t)\|_{\ell^2}\le \|u(t)\|_{\infty}$ by the discussion in
Section~\ref{sec:ap}. Hence the $\ell^2$ norm remains bounded by Theorem~\ref{thm:app}.
Therefore, it remains to prove that the frequency module $\cL^*$ is preserved.

 We start with the quasiperiodic case where $\cL^*$ is finitely-generated. Then, it follows from Proposition~\ref{characterisation}) that the frequency module is preserved. 
Indeed, if $u_0 \in AP(\cL^*)$ then, for every sequence $h_j$ with $Ah_j\to 0$, we have $u_0( h_j + \cdot)-u_0( \cdot) {\to 0}$ uniformly. But then by \eqref{dif}, $u(t, h_j + \cdot)-u(t, \cdot) \to 0$ uniformly as well and thus $u(t) \in AP(\cL^*)$ 
by Proposition~\ref{characterisation}. 
The general case is reduced to this particular one since any function $u\in AP(\cL^*)$ can be approximated in the uniform topology by quasiperiodic functions $u_n\in AP(\cL_n^*)$, where $\cL_n^*\subset \cL^*$ are finitely-generated.

The proof that spatial derivatives are almost-periodic follows analogously, the only difference is to use~\eqref{smo} instead of \eqref{dif}.  Temporal derivatives can then be expressed through the spatial ones using equation~\eqref{eq:SH}.
 \end{proof}

\begin{rmk} \label{rmk:globalAP} 
	Let $H\subset\OO(d)$ be a holohedry.
	By $\E(d)$-equivariance, $H$-invariance of solutions is preserved.
Hence it follows from Proposition~\ref{prop:globalAP} that if $u_0\in AP(\cL^*_H)$, then 
$u(t)\in AP(\cL^*_H)$ for all $t>0$.
\end{rmk}

\subsection{Global existence in $\ell^1(\cL^*)$} \label{sec:globall1}

Let $\cL^*\subset\R^d$ be a countable subgroup  and
suppose that $\cL^*$ is not uniformly discrete.
Proving  global existence of almost-periodic solutions in $\ell^1(\cL^*)$   is much more delicate than in $AP(\cL^*)$; in particular, there are no general methods to obtain globally bounded/dissipative solutions in $\ell^1(\cL^*)$ even when $\cL^*$ is finitely generated. 
Moreover, even for simple examples like the Swift-Hohenberg equation, there is no known result that excludes the possibility that typical solutions in $\ell^1(\cL^*)$ blow up in finite time.

However, when $\cL^*$ is finitely generated, there are methods for proving that quasiperiodic solutions with sufficient initial smoothness remain defined in $\ell^1(\cL^*)$ though without a bound on the $\ell^1$ norm (we can only check that the $\ell^1$ norm grows at most exponentially).
For equation~\eqref{eq:SH} 
we adapt  below one of such approaches, see also  \cite{Giga05, GigaMahalovYoneda11,M98,M99} and  references therein.
In particular, our results apply to quasicrystal solutions in $\ell^1(\cL^*_H)$.

Suppose that $\cL^*$ is finitely generated and that $u_0$ is an initial condition in $\ell^1(\cL^*)$. Then $u_0\in AP(\cL^*) \subset \ell^2(\cL^*)$ and the global well-posedness in $C_b(\R^d)$ implies that $u(t)\in AP(\cL^*) \subset \ell^2(\cL^*)$ for all $t\geq 0$ by Proposition \ref{prop:globalAP}.

To avoid blow up in finite time, we make some additional assumptions on the decay of the amplitudes of the initial condition $u_0$. These conditions are defined in terms of the hull function $U_0$ of $u_0$. Recall that this hull function is an element of $C(\T^p) \subset L^2(\T^p)$ where $p$ is the number of $\Z$-independent vectors that generate $\cL^*$.
For $U\in L^2(\T^p)$ we define
$$
\|U\|_{H^s} :=
\Big( \sum_{m\in\Z^p} ( |m|^2+1 )^s |U_m|^2 \Big)^{\frac{1}{2}}
 $$
to be its Sobolev norm of degree $s$, and we denote by $H^s(\T^p)$ the Hilbert space for which this norm is finite. 

\begin{prop} \label{prop:Hs}
Let $u\in AP(\cL^*)$ with hull function $U\in C(\T^p)$.
Suppose that $U\in H^s(\T^p)$ for some $s>p/2$.
Then $u\in \ell^1(\cL^*)$ and there is a constant $C>0$ independent of $u$ such that
$\|u\|_{\ell^1}\le C\|U\|_{H^s}$.
\end{prop}

\begin{proof}
By Cauchy-Schwarz,
$$
\|u\|_{\ell^1}  =\sum_{m\in\Z^p}|a_m|
 \le \bigg(\sum_{m\in\Z^p}\frac1{(|m|^2+1)^s}\bigg)^{1/2}
\bigg(\sum_{m\in\Z^p}(|m|^2+1)^s |a_m|^2\bigg)^{1/2}\le
C\|U\|_{H^s}
$$
for $s>p/2$.
\end{proof}

The next step is to rewrite the Swift-Hohenberg equation for $u$ in terms of the hull function $U$. 
Recalling that $u(x) = U(Ax)$, it follows that  if $u$ and $U$ are sufficiently smooth, then
\begin{equation}\label{eq:tildedelta}
	\Delta u(x) =\bigg(\sum_{i=1}^d \Big( \sum_{j=1}^p A_{ji}  \frac{\partial}{\partial \phi_j}\Big) ^2 U\bigg)  (Ax):=(\widetilde \Delta U)(Ax) .
\end{equation}
Thus, $u$ satisfies the Swift-Hohenberg equation if and only if  $U$ satisfies
\begin{align} \label{eqU}
U_t = -(\widetilde{\Delta}+1)^2U+\lambda U-U^3. 
\end{align}
It is immediate that $\|U(t)\|_{L^2} \leq \|u(t)\|_{\infty}$ for all $t$, so solutions of \eqref{eqU} remain in $L^2(\T^p)$ for all time. 
In contrast to $\Delta$, the operator $\widetilde \Delta$ is not uniformly elliptic, but highly degenerate, so we cannot expect strong smoothing properties for the solutions of \eqref{eqU}. However, when the initial condition is smooth, the smoothness is preserved for all time:

\begin{thm}\label{thm:Hs} 
Suppose that the initial condition $U_0$ to equation \eqref{eqU} satisfies $U_0\in H^s(\T^p)$ for some $s\ge1$. Then the solution $U(t)$  belongs to $H^s(\T^p)$ for all $t\ge0$
and $\|U(t)\|_{H^s}$ grows at most exponentially in time.

Consequently, if $U_0\in H^s(\T^p)$ for some $s>p/2$, then $\|u(t)\|_{\ell^1}$ grows at most exponentially in time.
\end{thm}

\begin{proof} 
The second statement is immediate from the first by Proposition~\ref{prop:Hs}.
We give below only the formal derivation for estimates of $U$ in the $H^s$ norms; these can be justified in a standard way using e.g.\ Galerkin approximations. Also, we give these derivations only for $s=1$ and $s=2$; the general case is handled similarly.

Let us start with $s=1$ and $U_0\in H^1(\T^p)$ and let $U^i=\partial_{\phi_i}U$. Then $U^i$ satisfies
$$
\partial_tU^i=-(\widetilde \Delta+1)^2U^i+\lambda U^i-3U^2U^i,\quad  \ U^i|_{t=0}=U^i_0.
$$
Multiplying this equation by $U^i$, integrating over $\T^p$ and using that the operator $(\widetilde \Delta+1)^2$ is nonnegative, we arrive at
$$
\frac12\frac d{dt}\|U^i(t)\|^2_{L^2}\le\lambda\|U^i(t)\|^2_{L^2}
$$
and therefore $\|U^i(t)\|_{L^2}^2\le e^{2\lambda t}\|U^i(0)\|^2_{L^2}$. 
Thus, $\|U(t)\|_{H^1}$ grows at most exponentially in time.

Suppose now that $s=2$ and  $U_0\in H^2(\T^p)$ and denote $U^{ij}:=\partial_{\phi_1}\partial_{\phi_j}U$. Then
$$
\partial_tU^{ij}=-(\widetilde \Delta +1)^2U^{ij}+\lambda U^{ij}-3U^2U^{ij}-6UU^iU^j,\, \ \ U^{ij}|_{t=0}=U^{ij}_0.
$$
Again, the strategy is to multiply this equation by $U^{ij}$ and integrate over $\T^d$. The only difficult term to estimate is the last term in the right-hand side. 
The interpolation inequality $\|U\|_{W^{1,4}}\le C\|U\|_{\infty}^{1/2}\|U\|^{1/2}_{H^2}$ gives
$$
\Big| \int_{\T^p}UU^iU^jU^{ij}\,d\phi \Big|\le \|U\|_{\infty}\|\nabla_\phi U\|_{L^4}^2\|U\|_{H^2}\le C\|U\|_{\infty}^2\|U\|^2_{H^2}\le C_1\|U\|^2_{H^2}
$$
where in the last inequality we used that $\|U\|_{\infty} = \|u\|_{\infty}$ is uniformly bounded in time. Multiplying the  equation by $U^{ij}$ and summing over $i$ and $j$, 
$$
\frac12\frac d{dt}\|U(t)\|^2_{H^2}\le (C+\lambda)\|U(t)\|^2_{H^2}\, .
$$
Again we conclude that $\|U(t)\|_{H^2}$ grows at most exponentially in time.
\end{proof}

\section{The Swift-Hohenberg equation; model-specific results}\label{sec:SH2}

In this section, we analyse the Swift-Hohenberg equation~\eqref{eq:SH} in more detail. In Subsection~\ref{sec:A}, we discuss further properties of the global attractor $\cA_\lambda$. 
In particular, for $d=2$ we show that $\cA_0=\{0\}$ at the bifurcation point $\lambda=0$, see Theorem~\ref{thm:A}, which seems to be new. Our proof of this result is restricted to the planar case although we expect that it remains true for higher dimensions as well.

In Subsection \ref{sec:grad}, for $d\le 9$, we construct in Theorem~\ref{thm:SH} a family $u_\lambda(t)$ of time-varying  quasicrystal solutions  of equation~\eqref{eq:SH} for $t\ge0$ with $\ell^2$ norm growing like $\sqrt\lambda$. In addition, for $\lambda\in(0,1]$, we prove that these solutions remain bounded away from the set of spatially constant functions.
To obtain these results, we utilise the gradient structure of equation~\eqref{eq:SH} restricted to the set of spatially almost-periodic functions. To verify property~(iii) of quasicrystal, we use the analyticity in space and time of solutions of equation~\eqref{eq:SH}. 

In Subsection~\ref{sec:attr}, for $d\le 9$,
we construct a different family of solutions $u_\lambda(t)$ which are defined for all $t\in\R$ (and therefore belong to $\cA_\lambda$) and tend to zero exponentially as $t\to-\infty$. We do not know how to verify property~(iii) for such solutions.

\subsection{The global attractor $\cA_\lambda$}
\label{sec:A}

The key result of this subsection is the triviality of the global attractor $\cA_\lambda$ at $\lambda=0$ when $d=1$ and $d=2$.

\begin{thm}\label{thm:A} For all $d\ge1$ and $\lambda<0$, the attractor $\cA_\lambda$ of the Swift-Hohenberg equation~\eqref{eq:SH} is trivial: $\cA_\lambda=\{0\}$, and $\cA_\lambda$ is exponentially attracting in $L^2_b(\R^d)$.

Moreover, $\cA_0=\{0\}$ for $d=1$ and $d=2$.
\end{thm}

 \begin{proof}
As in Appendix~\ref{app:global}, 
 we consider the weight function $\varphi(x)=\varphi_{\eps,x_0}(x)=\exp\{-\sqrt{\eps^2|x-x_0|^2+1}\}$. 
Denote the $L^2$-inner product by $(\;,\;)$.
For $\eps$ sufficiently small, it follows from~\eqref{eq:L2}
and Lemma~\ref{lem:lin} that solutions to~\eqref{eq:SH} satisfy
$$
\frac d{dt}(u^2,\varphi)\le (C\eps^2+2\lambda)(u^2,\varphi)-2(u^4,\varphi). 
$$
In particular, for $\lambda<0$ we can shrink $\eps$ so that
$\frac d{dt}(u^2,\varphi)\le \lambda(u^2,\varphi)$.
It follows that
$(u^2(t),\varphi)\le e^{\lambda t}(u_0^2,\varphi)$ for $t\ge0$.
Taking the supremum over $x_0$, it follows from~\eqref{eq:L2b} that
$\|u(t)\|_{L^2_b}^2\le Ce^{\lambda t}\|u_0\|_{L^2_b}^2$ for $t\ge0$.
Hence $\cA_\lambda=\{0\}$ and is exponentially attracting.

Next, we specialise to $\lambda=0$ and $d\le2$.
For $\lambda=0$, equation~\eqref{eq:SH} reads
\begin{align} \label{eq:SH0}
u_t=F(u, 0)=-(\Delta+1)^2u-u^3,
\end{align}
and we have
$$
\frac d{dt}(u^2,\varphi)\le C\eps^2(u^2,\varphi)-(u^4,\varphi). 
$$
It follows from~\eqref{eq:u4} that 
$$
	\frac d{dt}(u^2,\varphi)+\eps^2(u^2,\varphi)\le (C+1)\eps^2(u^2,\varphi)-C'\eps^d(u^2,\varphi)^2.
	$$
Since $d\le2$, 
$$
	\frac d{dt}(u^2,\varphi)+\eps^2(u^2,\varphi)\le 
	 -C'\eps^2\{(u^2,\varphi)-C''\}^2+\eps^2C'''
	 \le \eps^2C'''.
	$$
Here, $C,C',C'',C'''>0$ are constants independent of $\eps$.
This inequality guarantees that for any solution $u(t)$ of \eqref{eq:SH0}, there exists $t_0=t_0(u,\eps)>0$ such that
\begin{equation}\label{estimate}
\int_{\R^d}|u(t)(x)|^2e^{-\sqrt{\eps^2|x|^2+1}}\,dx \le 2C'''
\quad\text{for all $t\ge t_0$}.
\end{equation}
In particular, any solution $u(t)$ in the attractor $\cA_0$ satisfies \eqref{estimate} for all $t\in\R$. Therefore, by the monotone convergence theorem, $\cA_0$ is a bounded set in $L^2(\R^d)$.

We now proceed to study equation~\eqref{eq:SH0} in $L^2(\R^d)$. The Swift-Hohenberg equation is globally well-posed in $L^2(\R^d)$, possesses the parabolic smoothing property, and the Lyapunov function
$
\cP_0(u)=\tfrac12\|(\Delta+1)u\|^2_{L^2}+\tfrac14\|u\|^4_{L^4} \geq 0
$
for which we have
\begin{equation} \label{eq:grad}
\frac d{dt}\cP_0(u(t))=-\|F(u(t),0)\|^2_{L^2}  =-\|u_t(t)\|^2_{L^2}
\end{equation}
along solutions to~\eqref{eq:SH0} in $L^2(\R^d)$.
Due to the parabolic smoothing property for $L^2$ solutions,
$$
\cP_0(u(t_0+1))\le C(\|u(t_0)\|_{L^2}^2+\|u(t_0)\|^4_{L^2}),
$$
so integrating~\eqref{eq:grad} we obtain
\begin{equation} \label{eq:diss}
\int_{t_0+1}^t\|u_t(s)\|^2_{L^2}\,ds = \cP_0(u(t_0+1)) - \cP_0(u(t)) \le  \cP_0(u(t_0+1))  \le C\|u(t_0)\|^2_{L^2}.
\end{equation}
Moreover, differentiating~\eqref{eq:SH0} with respect to $t$ and multiplying by $u_t$, we get that
$
\frac d{dt}\|u_t(t)\|_{L^2}^2\le 0.
$
Hence $t\to\|u_t(t)\|^2_{L^2}$ is decreasing, and the dissipation integral~\eqref{eq:diss} gives
$$
\|u_t(t)\|^2_{L^2}\le C\frac{\|u(t_0)\|_{L^2}^2}{t-t_0-1}.
$$
Using that $u(t)$, $t\in\R$, is on the attractor (which is bounded in $L^2$) and passing to the limit $t_0\to-\infty$, we conclude that $u_t(t)\equiv 0$, so $u$ is a steady-state solution.

Finally, setting $u_t=0$ in equation~\eqref{eq:SH0} and multiplying by $u$ gives that $\|u\|_{L^4}=0$. Hence the only steady-state solution is $u\equiv0$, so $\cA_0=\{0\}$.
\end{proof}

 For $\lambda>0$, it is well-known that the  attractor $\cA_\lambda$ is now highly nontrivial, containing uncountably many group orbits of spatially-periodic steady-state solutions (see e.g.\ \cite{CrossHoh, Satt79}). Moreover, $\cA_\lambda$ has a very rich spatially chaotic structure since it contains an infinite dimensional  essentially unstable manifold of the zero equilibrium, which is parametrised by the proper space of frequency modulated functions, see \cite{Zel04} for the general theory and applications to reaction-diffusion equations. 
However, provided $d=2$ (or $d=1$), the size of the attractor $\cA_\lambda$ remains small for small positive $\lambda$ and tends to zero as $\lambda\to0$.

\begin{cor} \label{cor:A} Suppose that $d=2$. Then
$\lim_{\lambda\to0^+} \|\cA_\lambda\|_{L^\infty}=0$.
\end{cor}

\begin{proof} This estimate is a standard corollary of the upper semicontinuity of attractors $\cA_\lambda$ at $\lambda=0$ in $L^2_{loc}(\R^2)$ and the spatial invariance of the equation. Indeed, since $\cA_0=\{0\}$, due to upper semicontinuity, for every $\eps>0$, there exists $\lambda_\eps>0$ such that
$$
\|\cA_\lambda-0\|_{L^\infty([-1,1]^2)}\le\eps,\ \ \lambda<\lambda_\eps.
$$
Since equation~\eqref{eq:SH}, and therefore the attractor $\cA_\lambda$, is translation-invariant the last formula implies that
$$
\|\cA_\lambda\|_{L^\infty(\R^2)}=\sup_{h\in\R^2}\|\cA_\lambda\|_{L^\infty([-1+h_1,1+h_1]\times[-1+h_2,1+h_2)}=\|\cA_\lambda\|_{L^\infty([-1,1]^2)}\le\eps,\ \ \lambda<\lambda_\eps
$$
and the corollary is proved.
\end{proof}

\begin{rmk}
We do not have a convergence rate for the limit in Corollary~\ref{cor:A}  nor a rate of attraction to zero at the bifurcation point $\lambda=0$. 
\end{rmk}

\subsection{Gradient structure}\label{sec:grad}
In this section, we study further the Swift-Hohenberg equation~\eqref{eq:SH} for $1\le d\le 9$ and show how to exploit the gradient structure.
 Recall that equation~\eqref{eq:SH} formally possesses a Lyapunov functional
$$
\cL (u):=\tfrac12\|(\Delta+1)u\|^2_{L^2}+\tfrac14\|u\|^4_{L^4}-\tfrac12\lambda\|u\|^2_{L^2},\ \ \frac d{dt}\cL(u(t))=-\|\partial_tu(t)\|^2_{L^2}\le0.
$$
This gives a genuine Lyapunov functional in a bounded domain or in the case of square integrable initial data, but in the case where the solution is only spatially bounded, these integrals may be infinite, so we do not have a true gradient structure. Such systems are usually referred 
as \emph{extended gradient systems}, see \cite{Gal01,Zelik03,Zel04} for more details. 
In particular, the gradient structure may be restored if we extend equation~\eqref{eq:SH} to the space $\cM_{sp}$ of spatially invariant Borel measures
on $C_b(\R^d)$ endowed by the topology of $C_{loc}(\R^d)$:
\begin{align}
\nonumber
& \cL (\mu):=\int\big(\tfrac12|(\Delta+1)u|^2+\tfrac14|u|^4-\tfrac12 \lambda|u|^2\big)\mu(du),\\ 
& \frac d{dt}\cL(\mu(t))=-\int|\partial_tu(t)|^2\mu(du)\le0,
\label{meas} 
\end{align}
see \cite{Zelik03,Zel04}. In the case of almost-periodic initial data $u_0\in AP(\cL^*)$, using the Lebesgue (Haar) measure on the hull of $u_0$, we may rewrite \eqref{meas} in the form 
\[
P_\lambda(u)=\tfrac12\|(\Delta+1)u\|_{\ell^2}^2-\tfrac12\lambda \|u\|_{\ell^2}^2
+\tfrac14 \|u^2\|_{\ell^2}^2\, .
\]
Indeed, this quantity exists provided $u, \Delta u\in AP(\cL^*)$. Hence, due to the smoothing property for equation~\eqref{eq:SH}, $P_\lambda(u(t))$ is well-defined for $t>0$ for all solutions $u(t)$ such that $u_0\in AP(\cL^*)$. 

Although the fact that $P_\lambda(u)$ is a Lyapunov functional follows directly from \eqref{meas}, we give below an independent proof.

\begin{prop} \label{prop:inner}   Let $\cL^*\subset\R^d$ be a countable subgroup. Then
\begin{itemize}
\item[(a)]
$\langle u,vw \rangle_{\ell^2} = \langle uv,w \rangle_{\ell^2}$
for all $u,v,w\in AP(\cL^*)$.
\item[(b)]
$\langle (\Delta+1)^2u,v\rangle_{\ell^2} =
\langle (\Delta+1)u,(\Delta+1)v\rangle_{\ell^2}$
for all $u,v\in AP(\cL^*)$ such that also $(\Delta+1)^2u,(\Delta+1)v\in AP(\cL^*)$.
\item[(c)] $\langle u,u\rangle^2_{\ell^2}\le \langle u^2,u^2\rangle_{\ell^2}$
for all $u\in AP(\cL^*)$.
\end{itemize}
\end{prop}

\begin{proof}  Using approximation arguments, it is enough to verify these formulas for finitely generated subgroups $\cL^*$. Recall from Section~\ref{sec:AP} that $\langle u,v \rangle_{\ell^2} = \langle U, V\rangle_{L^2}$ where $U,\, V:\T^p\to\R$ are the hull functions of $u$ and $v$.

\indent To prove~(a),  note that $u,v,w \in AP(\cL^*)$ implies that $uv, vw \in AP(\cL^*)$ and hence $\langle u,vw \rangle_{\ell^2} = \langle U, VW\rangle_{L^2} = \langle UV, W\rangle_{L^2} = \langle uv,w \rangle_{\ell^2}$.

To prove~(b), note that the hull function of $(1+\Delta)u$ is given by $(1+\widetilde \Delta)U$, with $U$ the hull function of $u$ and the operator $\widetilde \Delta$ as given in \eqref{eq:tildedelta}. As a result, $$\langle (\Delta+1)^2u,v\rangle_{\ell^2} = \langle  (\widetilde\Delta+1)^2U,V\rangle_{L^2} =  \langle  (\widetilde \Delta+1)U,(\widetilde \Delta + 1)V\rangle_{L^2} = \langle (\Delta+1)u,(\Delta+1)v\rangle_{\ell^2}\, ,$$
where the second equality follows from integration by parts on $\T^p$.

Finally, by Cauchy-Schwarz,
$$
\langle u,u\rangle^2_{\ell^2} = \langle U, U\rangle^2_{L^2} =  \langle U^2,1\rangle^2_{L^2} \le \langle U^2,U^2\rangle_{L^2}\langle1,1\rangle_{L^2} =\langle u^2,u^2\rangle_{\ell^2}. 
$$
yielding part~(c).
\end{proof}

\begin{prop}  \label{prop:P}
Let $u_{\lambda}(t)$ be a solution to the equation~\eqref{eq:SH} and assume that $u_{\lambda}(0) \in AP(\cL^*)$.
Then
\begin{equation}\label{AP-Lyap}
\frac{d}{dt} P_\lambda(u_{\lambda}(t))= - \|F(u_{\lambda}(t),\lambda)\|_{\ell^2}^2=-\|\partial_tu_\lambda(t)\|^2_{\ell^2}
\quad\text{for all $t>0$.}
\end{equation}
\end{prop}

\begin{proof}
Let $u_\lambda(0)\in AP(\cL^*)$ and $t>0$. 
Due to the smoothing property, $u(t), \Delta u(t), \Delta^2 u(t),\partial_tu(t),\Delta\partial_t u_t \in AP(\cL^*)$ and, by Proposition~\ref{prop:inner}{\it (a,b)}, we have
\begin{align*}
(dP_\lambda)_{u}\partial_t u & =\langle (\Delta+1)u,(\Delta+1)\partial_t u\rangle -\lambda\langle u,\partial_t u \rangle
+\langle u^2,u\partial_tu\rangle
\\ & =\langle (\Delta+1)^2u,\partial_tu\rangle -\lambda\langle u,\partial_t u \rangle
+\langle u^3,\partial_t u\rangle
=-\langle F(u,\lambda),\partial_t u\rangle.
\end{align*}
 Thus it follows as usual that
$$
\frac{d}{dt}P_\lambda(u_{\lambda}(t))=(dP_\lambda)_{u_{\lambda}(t)}\partial_tu_{\lambda}(t)=-\langle F(u_{\lambda}(t),\lambda),F(u_{\lambda}(t),\lambda)\rangle.
$$

\vspace{-5ex}
\end{proof}

\begin{rmk} \label{rmk:min}
When $\lambda\le0$, the potential function $P_{\lambda}$ has a unique global minimum at $u=0$.
Indeed,  by Proposition~\ref{prop:inner}{\it (c)},
$P_\lambda(u)\ge \frac14 \bigl(\|u\|_2^4-2 \lambda\|u\|_2^2 \bigr)$.
\end{rmk}

\begin{rmk}
 We note that \eqref{AP-Lyap} holds up to $t=0$ if $u_\lambda(0)$ is smooth enough. In particular, it is not difficult to prove that $t\mapsto P_\lambda(u_\lambda(t))$ is continuous at $t=0$ if $u_\lambda(0),\Delta u_\lambda(0)\in AP(\cL^*)$.
\end{rmk}

Our next results describe the behaviour of solutions of equation~\eqref{eq:SH} in the $\ell^2$ norm assuming that $u_0\in AP(\cL^*)$.
\begin{prop}\label{Prop-en} Let $u_0\in AP(\cL^*)$ and $N(u):=\|u(t)\|^2_{\ell^2}$. Then
\begin{equation}\label{enest}
\frac d{dt}N(u(t))\le N(u(t))(\lambda-N(u(t)),\ \ t>0.
\end{equation}
Moreover, the function $t\to N(u(t))$ is continuous at $t=0$.
\end{prop}

\begin{proof} 
Using Proposition~\ref{prop:inner}, we compute that
\begin{align*}
\frac12\frac d{dt} N & =\langle u,u_t\rangle = \langle u,-(\Delta+1)^2u\rangle
+\langle u,\lambda u\rangle + \langle u,-u^3 \rangle \\
& = -\langle (\Delta+1)u,(\Delta+1)u\rangle
+\lambda \langle u,u\rangle - \langle u^2,u^2 \rangle \\
& \le \lambda   \langle u,u\rangle -  \langle u,u\rangle^2 =N(\lambda - N).
\end{align*}

Continuity at $t=0$ follows from the fact that the function $t\to u(t)$ is continuous in the $C_b$ norm at $t=0$ for any $u_0\in AP(\cL^*)$. 
\end{proof}

\begin{cor}\label{Cor3.l2} Let $u_0\in AP(\cL^*)$ with corresponding solution $u(t)\in AP(\cL^*)$.
 
(a) For $\lambda<0$, we have exponential convergence to zero in the $\ell^2$ norm:
$$
\|u(t)\|_{\ell^2}\le e^{\lambda t}\|u_0\|_{\ell^2},\ \ t\ge0;
$$
 
(b) At the bifurcation point $\lambda=0$, we have polynomial convergence to zero:
$$
\|u(t)\|^2_{\ell^2}\le \frac{2\|u_0\|^2_{\ell^2}}{2\|u_0\|^2_{\ell^2}+1},\ \ t\ge0;
$$
 
(c) For $\lambda>0$, we have an absorbing ball of radius proportional to $\sqrt{\lambda}$ in $\ell^2$. Namely, for any $\eps>0$, there exists $T=T_\eps$ (independent of $u_0$) such that
\begin{itemize}
\item[]
\hfill \qquad\qquad $\|u(t)\|_{\ell^2}\le (1+\eps)\sqrt{\lambda},\ \ t\ge T.$ \hfill \qed
\end{itemize}
\end{cor}

\begin{thm} \label{thm:SH}
Suppose that $1\le d\le 9$.
Then for each holohedry $H\subset\OO(d)$, there is a family of spatially quasiperiodic solutions $u_\lambda(t) \in\ell^1(\cL^*_H)$ to the Swift-Hohenberg equation, defined for all $t\ge0$, and a constant $C_H\in(0,1]$ such that
\begin{equation} \label{eq:sqrt}
C_H\sqrt\lambda\le \|u_\lambda(t)\|_2\le \sqrt\lambda \quad\text{for all $t\ge0$, $\lambda>0$}.
\end{equation} 
These solutions are quasicrystals if $\cL^*_H$ is not uniformly discrete.

For $\lambda\in(0,1]$, these solutions are bounded away from being spatially constant:
\begin{equation} \label{eq:sep}
\inf_{t\ge0}\inf_{c\in\R}\|u_\lambda(t)-c\|_{L^\infty}  >0 .
\end{equation} 

Moreover, for each $\lambda>0$, the solution $u_\lambda(t)$ satisfies condition~(iii) from the introduction for  all but at most countably many $t>0$.
\end{thm}
 
\begin{proof}
As in Part~I, we consider 
initial conditions $u_0(\lambda)$ of the form
\[
u_\lambda(0)(x)=a\sqrt\lambda \sum_{\gamma\in H} e^{i\gamma k_0\cdot x},
\]
where $k_0$ is a fixed unit vector in $\R^d$ and $a>0$ is chosen
appropriately.

Since $u_\lambda(0)$ is a finite sum, it is clear that $u_\lambda(0)\in \ell^1(\cL^*_H)$ and that its hull function
$U_\lambda(0)$ is in $H^s(\T^p)$ for all $s\geq 0$. By Theorem~\ref{thm:Hs}, 
$u_\lambda(t)\in \ell^1(\cL^*_H)$ for all $t\ge0$.
 
Estimates~\eqref{eq:sqrt} and~\eqref{eq:sep} are proved in 
Part~I (Theorem~3.1(c,d)).
Thus, it remains to show that $u_\lambda(t)$
satisfies condition~(iii) for all but countably many values of~$t$. To this
end, we perturb slightly the initial data $u_\lambda(0)$ in such a way that
$a_{\lambda,k}(0){\neq0}$ for all $k\in \cL^*_H$ and so that the various estimates
remain essentially unchanged.
Since the solution $u_\lambda(t)$ is analytic in  $t>0$, we conclude
that the amplitudes $a_{\lambda,k}(t)$ are analytic. Since, in addition, they are
continuous at $t=0$, we get that $a_{\lambda,k}(t)$ is not identically zero
for each $k\in \cL^*_H$, so $a_{\lambda,k}(t)=0$ for some $k$ 
for at most countably many $t>0$. This completes the proof of the theorem.
\end{proof}

\begin{rmk} The condition $\lambda\le1$ is essential for the separation~\eqref{eq:sep}. Indeed, if $\lambda>1$, then equation~\eqref{eq:SH} possesses two further spatially constant equilibria $u(t)\equiv \pm\sqrt{\lambda-1}$, which are exponentially stable, so \emph{a priori} it could be that the solution $u_\lambda(t)$ converges to one of them as $t\to\infty$.
\end{rmk}

\subsection{Quasicrystals in $\cA_\lambda$}
\label{sec:attr}

Let $\lambda>0$. Corresponding to each holohedry $H$,
it is possible to construct spatially quasiperiodic solutions $u_\lambda(t)\in\cA_\lambda$, defined for all $t\in\R$, which converge to zero as $t\to-\infty$. Indeed, for $\lambda>0$, we may construct an essentially unstable manifold $\cM_\lambda\subset C_b(\R^d)$, which consists of solutions approaching zero sufficiently fast as $t\to-\infty$ and which is parametrised by the solutions of the linear problem
$$
\partial_t v+(\Delta+1)^2u=\lambda v, \ v|_{t=0}=v_0, \ t\le0,\ \ v(t)\to0 \ \ \text{sufficiently fast}.
$$
Namely, the solution on this local manifold is defined via
$$
u(t)=v(t)+M_\lambda(v(t)),
$$
where the function $M_\lambda:C_b(\R^d)\to C_b(\R^d)$ is smooth and satisfies $M_\lambda(0)=M'_\lambda(0)=0$, see~\cite{Zel04} for the details.
 
Moreover, if we take $v_0=u_\lambda(0)$, where $u_\lambda(0)$ is as constructed in Theorem~\ref{thm:SH}, then the solution $v(t)$ is almost-periodic for all $t\le0$ ($v(t)\in AP(\cL^*_H)$), approaching zero sufficiently fast and  satisfying $P_\lambda(v(t))<0$ for all small $t$. From this we may conclude that $u_\lambda(t)\in AP(\cL^*_H)$ for all $t\in\R$ and $P_\lambda(u_\lambda(t))<0$ for all $t$.
 
In contrast to the solutions obtained in Theorem~\ref{thm:SH}, we cannot freely perturb \emph{all} amplitudes $a_{\lambda,k}(0)$ since we need to be sure that $v(t)$ belongs to the base of the essentially unstable manifold. Hence, it is not clear whether property~(iii) from our definition of quasicrystal is satisfied although we argue that the failure of this property is a pathology and that it would hold for typical examples in the highly unlikely event that it failed for the specific equation~\eqref{eq:SH}.
 
Indeed,  we can write $u_\lambda(t)(x)=\sum_{k\in\cL^*_H}a_{\lambda,k}(t)e^{ik\cdot x}$ where
\[
\frac{d}{dt} a_{\lambda,k}(t)
=(\lambda-(|k|^2-1)^2)a_k(t)-\sum_{k'+k''+k'''=k}a_{\lambda,k'}(t)a_{\lambda,k''}(t)a_{\lambda,k'''}(t).
\]
Since the amplitudes $a_{\lambda,k}$ are analytic, for property~(iii) to fail we would need
$$
\sum_{k'+k''+k'''=k}a_{\lambda,k'}(t)a_{\lambda,k''}(t)a_{\lambda,k'''}(t)\equiv0
$$
for ``many'' $k\in\cL^*_H$ which seems unlikely since we know that $a_{\lambda,k}$ do not vanish identically for the frequencies $k$ such that $|k|$ are close to one and the system of cubic equations written above somehow couples all modes. Unfortunately, we do not have a rigorous proof for property~(iii); we plan to return to this problem somewhere else.

\section{The Brusselator model}
\label{sec:beyond}

Finite wavelength instabilities frequently occur for many other interesting equations, for instance in reaction diffusion systems where these are usually referred to as Turing instabilities. A prominent example is the so-called Brusselator, which is the main object for consideration in this section. It is given by the system of equations
\begin{equation}\label{eq:Brusselator}
\begin{cases}
\partial_tu = d_1 \Delta u + A - (B+1)u+ u^2 v,\ \ u|_{t=0}=u_0\\
\partial_t v = d_2\Delta v +Bu - u^2v,\ \ v|_{t=0}=v_0,
\end{cases}
\end{equation}
 where $x\in \R^d$ and $d_1,d_2,A,B$ are positive parameters. 
The system has been studied broadly from a physical viewpoint, in particular deriving conditions for instability and bifurcations of simple patterns, e.g. \cite{Verdasca}.

Briefly, it can be readily computed that the spatially homogeneous steady-state $(\bar u,\bar v) = (A,B/A)$ undergoes a generic Turing instability at $B=B_c:=(1+A \eta)^2$ provided that $1<\eta:= \sqrt{d_1/d_2}<((1+A^2)-1)/A$. 
Then the linearised equation for this steady-state at the bifurcation point  lies in the closed left half plane and has zero real part only at the origin. Moreover, the critical eigenfunctions are of the form $(u_0,v_0)e^{i k\cdot x}$ for  wavevectors $k\in\R^d$ of length $k_c = \sqrt{A/\eta}$, and where $(u_0,v_0)$ is a multiple of
$(-A,A\eta^2+\eta)$.
The critical spectrum for $B\approx B_c$ is then analogous to that of the Swift-Hohenberg equation~\eqref{eq:SH} for $\lambda\approx 0$.

Local well-posedness of classical solutions in various spaces including $C_b(\R^d)$ is standard for this semi-linear parabolic system. 
Note that equation~\eqref{eq:Brusselator} preserves the cone $\{u\ge0,v\ge0\}$ 
due to the form of reaction kinetics (the non-diffusive part of the right-hand side),  see e.g.\  \cite{Hollis}. 
These are the physically relevant solutions and hence we restrict attention to this cone.

There are two essential differences in comparison to equation~\eqref{eq:SH}.
First, the steady-state $(\bar u,\bar v)$ is not globally asymptotically stable before the bifurcation point (for $B\approx B_c$, $B<B_c$). Hence, the global attractor $\cA_B$ is nontrivial even for $B<B_c$ and only local stability can be proved.
Second, there is no gradient-like structure
  and we are unable to check  the separation of a quasicrystal solution from spatially constant solutions (the analogue of \eqref{eq:sep} is unclear here). \emph{A priori}, such solutions may be homoclinic to $(\bar u,\bar v)$ in $L^\infty(\R^d)$.

Despite these differences, many of our  key results for equation~\eqref{eq:SH} can be translated to the case of Brusselator model~\eqref{eq:Brusselator}. 
In particular, we prove existence of a global attractor in $L^\infty(\R^d)$, extending the arguments of~\cite{Rothe} to the case of unbounded domains and infinite energy solutions. This result is true for all admissible values of the parameters and any space dimension $d$. 
Also, we can find quasicrystal solutions defined for all $t\ge0$ and satisfying condition (iii), as well as quasicrystal solutions defined for all $t\in\R$ and converging to $(\bar u,\bar v)$ as $t\to-\infty$.

\subsection{Global existence}

Define $L^\infty_+(\R^d)$ to be the space of $L^\infty$ functions
$(u,v):\R^d\to [0,\infty)^2$.
 We start by verifying the global existence and dissipativity of the solution semigroup 
$$
S(t):L^\infty_+(\R^d)\to L^\infty_+(\R^d), \qquad  S(t)(u_0,v_0):=(u(t),v(t)).
$$
We note that previously, 
global existence on \emph{bounded} domains with various boundary conditions, and for a class of systems including the Brusselator, had been proven in, e.g., \cite{Hollis,Rothe}. In addition, the case $d_1=d_2$ (which does not admit a Turing instability) had been studied on $\R^d$ in \cite{Guo}.

\begin{thm}\label{Th4.main} For any $d\ge1$ and any $A,B,d_1,d_2>0$, the Brusselator system \eqref{eq:Brusselator} possesses a unique global solution $(u(t),v(t))\in L^\infty_+(\R^d)$ and the associated solution semigroup is locally Lipschitz continuous and is globally dissipative:
there exist $\alpha>0$, $C\in\R$ and a monotone increasing function $Q$ such that
\begin{equation}\label{B-dis}
\|(u(t),v(t))\|_\infty\le Q(\|(u_0,v_0)\|_\infty)e^{-\alpha t}+C.
\end{equation}
Here, $\alpha,\,C,\,Q$
are independent of $t$ and $(u_0,v_0)$, but may depend on $A,B,d_1,d_2$.
\end{thm}

\begin{proof}\textbf{for $d\le 3$.}
Global existence in $L^\infty_+(\R^d)$ can be shown by adjusting the $L^2$ estimates derived in \cite{Rothe} on bounded domains with Neumann boundary conditions to weighted spaces $L^2_\varphi(\R^d)$ with weight function $\varphi$ as described in the appendix. 
We sketch below the derivation of the basic dissipative estimate \eqref{B-dis}
 
{\it Step 1. Lower bounds for the $u$-component.} We apply the comparison principle to the first equation of \eqref{eq:Brusselator} dropping out the non-negative term $u^2v$. Solving explicitly the remaining linear heat equation with zero initial data, we end up with
\begin{equation}\label{lower1}
u(t,x)\ge \frac{A(1-e^{(B+1)t)})}{B+1}.
\end{equation}
Importantly, the solution $u$ becomes strictly positive and uniformly separated from zero for all $t>0$.
 
{\it Step 2. Upper bounds for the $v$-component.} Multiplying the second equation of~\eqref{eq:Brusselator} by $v$, we get the equation for $w=v^2$:
\begin{equation}\label{v2}
\partial_t w-\Delta w=-u^2v^2-(uv-B)^2-2|\nabla v|^2+B^2:=f(t,u,v).
\end{equation}
We estimate the function $f$ using \eqref{lower1} as follows:
$$
f(t,u,v)\le \begin{cases} B^2, &  t\le \frac{\ln2}{B+1}\\
                                      B^2-\frac{4A^2}{(B+1)^2}w, &  t\ge \frac{\ln2}{B+1}
                 \end{cases}
$$
and apply the comparison principle to \eqref{v2} to get
\begin{equation}\label{upper-v}
v^2(t,x)\le C (e^{-\alpha t}\|v_0\|^2_{\infty}+1),
\end{equation}
where the positive constants $C$ and $\alpha$ are independent of $t$, $u_0$ and $v_0$.
 
{\it Step 3. Weighted $L^2$ estimate for the $H^1$ norm of $v$-component.} To get this estimate, we multiply    equation \eqref{v2} by $\varphi$  and integrate with respect to $x\in\R^d$ and $t\in[T,T+1]$. 
Then, using~\eqref{upper-v} and~\eqref{eq:phi}, we end up with the desired estimate
\begin{equation}\label{grad}
\int_T^{T+1}(|\nabla v(t)|^2,\varphi)\,dt\le C_\eps(e^{-\alpha T}\|v_0\|^2_{\infty}+1),
\end{equation}
where $\eps>0$ is small enough and $C_\eps,\alpha>0$ are some positive constants.
 
{\it Step 4. Weighted $L^2$ estimates for the $u$-component.}
 Let us sum the equations for $u$ and $v$ in~\eqref{eq:Brusselator} to get
\begin{equation}\label{sum}
\partial_t (u+v) - \Delta(au+bv) = A-u.
\end{equation}
We multiply this equation by $(u+v)\varphi$, integrate by parts and use that
\begin{align*}
 & \big((d_1\nabla u+d_2\nabla v)  .(\nabla u+\nabla v),\varphi\big)
 \ge d_1(|\nabla u|^2,\varphi)+d_2(|\nabla v|^2,\varphi)-2(|\nabla u|.|\nabla v|,\varphi)
\\ & \qquad \qquad \qquad  \ge \tfrac{d_1}2(|\nabla u|^2,\varphi)+\Big(\big\{(\tfrac{d_1}2)^{1/2}|\nabla u|-\tfrac{d_1+d_2}{2d_1}|\nabla v|\big\}^2,\varphi\Big)-\tfrac{(d_1+d_2)^2}{2d_1}(|\nabla v|^2,\varphi)
\end{align*}
as well as
$$
(A-u,(u+v)\varphi)\le -\tfrac12(|u+v|^2,\varphi)+C(|v|^2,\varphi)+C_\eps
$$
which is valid for some positive $C$ and $C_\eps$. In this way, we obtain
$$
\frac d{dt}(|u+v|^2,\varphi)+\beta(|\nabla u|^2,\varphi)+\alpha(|u+v|^2,\varphi)\le C(|v|^2+|\nabla v|^2,\varphi)+C_\eps
$$
for some positive constants $\alpha,\beta,C_\eps$ and  $\eps>0$ sufficiently small (we estimate the terms containing the gradient of $\varphi$ exactly as in the appendix). Applying Gronwall's inequality to this relation and using \eqref{upper-v} and \eqref{grad}, we obtain
\begin{equation}\label{L2-u}
(|u(T)|^2,\varphi)+\int_{T}^{T+1}(|\nabla u(t)|^2,\varphi)\,dt\le C_\eps\left(e^{-\alpha T}(\|u_0\|_{\infty}^2+\|v_0\|_{\infty}^2)+1\right),
\end{equation}
for some positive constants $\alpha$ and $C_\eps$.
 
{\it Step 5. $L^\infty$ estimates for the $u$-component: Moser iterations.} 
Using the shifted weights $\varphi_{x_0}(x):=\varphi(x-x_0)$ in \eqref{L2-u} and taking the supremum over $x_0$, we derive uniformly local $L^2$ estimates for the $u$-component. Then, we analyse the equation for the $u$-component, namely, since the $v$-component is bounded in $L^\infty$, the nonlinearity $f(u,v):=u^2v$ actually has only quadratic growth rate (the presence of $v$ is not essential) and the critical exponent for the semilinear heat equation (where the $L^2$ energy estimates are known) is $p_c=1+4/d$. 

Therefore, for $d\le 3$, the quadratic nonlinearity is subcritical and the bootstrapping technique works. For instance, we may use Moser iterations in order to get the desired $L^\infty$ estimates. These iterations are very standard although a bit technical, so we omit the details, see also \cite{Hollis,Rothe,Guo}.
\end{proof}

\begin{rmk} The restriction $d\le3$ is necessary only for the bootstrapping in the last step of the proof.
This bootstrapping argument can be refined for $d=4$, but does not work for $d>4$. Nevertheless, the theorem remains true for all $d\ge1$. To prove this, we apply the mass dissipation condition \eqref{sum} following \cite{Kanel}, see also \cite{fell} where the global solvability for reaction-diffusion systems satisfying the mass dissipation condition is established for quadratic nonlinearities in any space dimension.
\end{rmk}

\begin{cor} Under the assumptions of Theorem \ref{Th4.main}, the solution $(u(t),v(t))$ becomes analytic in space and time in a strip for $t>0$, and the corresponding solution semigroup $S(t)$ possesses a locally compact global attractor $\cA$ which consists of analytic bounded functions in the corresponding strip.
\end{cor}

Indeed, the analytic smoothing property follows from the general results of \cite{TakTi} exactly as in the case of equation~\eqref{eq:SH} and the existence of a global attractor is an immediate corollary of the dissipative estimate \eqref{B-dis} and the parabolic smoothing property, see also \cite{Guo} where the existence of an attractor is established in the uniformly local spaces.

\subsection{Quasicrystals}
 
\begin{cor} If the initial data $(u_0,v_0)$ lies in $AP(\cL^*)$ for some countable subgroup $\cL^*\subset\R^d$, then the solution $(u(t),v(t))$ remains in
$AP(\cL^*)$ for all $t>0$ and is analytic in space and time.

If $H\subset\OO(d)$ is a holohedry and $\cL^*_H$ is not uniformly discrete, 
then there exist quasicrystal solutions $(u(t),v(t))$, $t>0$, in $[\ell^1(\cL^*)]^2$.
\end{cor}
These statements can be proved exactly as for equation~\eqref{eq:SH}.

\begin{rmk}
We see that most of the results obtained for equation~\eqref{eq:SH} have straightforward generalisations to the Brusselator system. In particular, for $B\approx B_c$, $B>B_c$, the spatially homogeneous steady-state $(\bar u,\bar v)$ becomes unstable and we have an essential unstable manifold. As in Section~\ref{sec:attr}, we may construct quasicrystals $u(t)$, $t\in\R$, on the global attractor $\cA$ that converge to $(\bar u,\bar v)$  as $t\to-\infty$. 

However, in the absence of a global Lyapunov function we are unable to prove existence of quasicrystals that satisfy the separation property~\eqref{eq:sep} from spatially constant solutions.
\end{rmk}

\appendix

\renewcommand{\thesubsection}{\Alph{section}.\arabic{subsection}}

\section{Global existence for the Swift-Hohenberg equation}
\label{app:global}

In this appendix, we consider the
global well-posedness and dissipativity of the Swift-Hohenberg equation~\eqref{eq:SH} in $L^\infty(\R^d)$ for moderate values of $d$.

\begin{thm} \label{thm:app}
Let $1\le d\le 9$, $\lambda_0>0$. There is a constant $C>0$ and a monotonic increasing function $Q$ such that
solutions of the PDE~\eqref{eq:SH} satisfy
$$
\|u(t)\|_{\infty}\le e^{-t}Q(\|u_0\|_{\infty})+C \quad 
\text{for all $t\ge0$, $\lambda\in(-\infty,\lambda_0]$}.
$$
\end{thm}

The proof of Theorem~\ref{thm:app} is a standard exercise for experts and uses well-trodden techniques, but it is not easy to pin down the precise details in the literature.
(For instance, global well-posedness in $L^2_b(\R^d)$ is claimed in~\cite{efendiev} for all $d\ge1$, with the main steps indicated very briefly for $1\le d\le 3$.) Since such techniques may not be so familiar for most readers of Part~I, we give complete details in $L^\infty(\R^d)$ for $d\le 7$ and indicate the additional steps for $d\le 9$. We do not know if Theorem~\ref{thm:app} holds for $d\ge10$.
(The result holds in $L^2_b(\R^d)$ for all $d\ge1$, but the smoothing argument from $L^2_b(\R^d)$ to $L^\infty(\R^d)$ in Subsection~\ref{sec:8} seems unclear for $d\ge10$.)

Unless stated otherwise, from now on, $C$ is a constant that may vary from line to line but which depends only on $d$ and $\lambda_0$.

\subsection{Estimates in weighted spaces}

For $\eps>0$, $x_0\in \R^d$, define
 the weight function
\[
\varphi_{\eps,x_0}:\R^d\to[0,1],
\qquad
\varphi_{\eps,x_0}(x)=\exp\{-\sqrt{\eps^2|x-x_0|^2+1}\}.
\]
For each $\varphi=\varphi_{\eps,x_0}$ we define the weighted $L^p$ spaces
$L^p_\varphi(\R^d)$ with norm
$\|u\|_{L^p_\varphi}=\|u\varphi^{1/p}\|_{L^p}$.
Similarly, we define the weighted Sobolev spaces
$H^k_\varphi(\R^d)$ with norm $\|u\|_{H^k_\varphi}=\|u\varphi^{1/2}\|_{H^k}$.
For each $d,\,n\ge1$, there is a constant $C>0$ such that
\begin{equation} \label{eq:phi}
| \nabla\varphi_{\eps,x_0}|\le C\eps\varphi_{\eps,x_0}, \quad |D^n_x\varphi_{\eps,x_0}|\le C\eps^n\varphi_{\eps,x_0}, \quad\text{for all $\eps>0$, $x_0\in\R^d$}.
\end{equation}
For $\eps>0$, $x_0\in\R^d$, let $\varphi=\varphi_{\eps,x_0}$.
Denote the $L^2$-inner product by $(\;,\;)$.

We multiply equation~\eqref{eq:SH} by $u\varphi$ and integrate over $\R^d$. 
Integrating by parts, this yields
\begin{equation} \label{eq:L2}
\frac12\frac{d}{dt}(u^2, \varphi) = -
((\Delta+1)u, (\Delta+1)(u\varphi))
+ \lambda(u^2,\varphi) -  (u^4, \varphi) .
\end{equation}

\begin{lemma} \label{lem:lin}
	There is a constant $C>0$ such that
\[
((\Delta+1)u, (\Delta+1)(u\varphi)) \ge
\tfrac12(|(\Delta +1)u|^2,\varphi)-C\eps^2(u^2,\varphi).
\]
for all $\eps>0$ sufficiently small.
\end{lemma}

\begin{proof}
Note that
$(\Delta+1)(u\varphi)=((\Delta+1) u)\varphi+u\Delta\varphi+2\nabla u\cdot\nabla\varphi$,
so
\[
((\Delta+1)u,  (\Delta+1)(u\varphi))=
 (|\Delta +1)u|^2,\varphi)
+A_1+A_2+2A_3
\]
 where
\[
  A_1   = (u,u\Delta\varphi), \quad
  A_2   = (\Delta u,u\Delta\varphi), \quad
  A_3= ((\Delta+1) u,\nabla u\cdot \nabla\varphi).
\]
By~\eqref{eq:phi}, $|A_1|= |(u^2,\Delta\varphi)|\le
C\eps^2 (u^2,\varphi)$.
Integration by parts yields
\begin{align*}
A_2
& = -\int \nabla u\cdot\nabla(u\Delta\varphi)
 = -\int |\nabla u|^2\Delta\varphi
-\int u\nabla u\cdot \nabla\Delta\varphi
\\ & = -\int |\nabla u|^2\Delta\varphi
-\tfrac12\int \nabla(u^2)\cdot \nabla\Delta\varphi
 = -\int |\nabla u|^2\Delta\varphi
+\tfrac12\int u^2\Delta^2\varphi.
\end{align*}
Hence by~\eqref{eq:phi},
 $|A_2|\le C\eps^2( |\nabla u|^2,\varphi)+C\eps^4 ( u^2,\varphi)$.

By the rough inequality $\sqrt{xy}\le x+y$ which holds for $x,y\ge0$,
\begin{align*}
|(\Delta+1)u\,\nabla u\cdot\nabla\varphi|
& \le \{|(\Delta+1)u|\,|\nabla\varphi|^{1/2}\}
\{|\nabla u|\,|\nabla\varphi|^{1/2}\}
\\ & 
\le |(\Delta+1)u|^2|\nabla\varphi|+
|\nabla u|^2|\nabla\varphi|.
\end{align*}
Hence it follows from~\eqref{eq:phi} that
\[
|A_3|\le C\eps (|(\Delta+1)u|^2,\varphi)+C\eps(|\nabla u|^2,\varphi)
\le \tfrac14 (|(\Delta+1)u|^2,\varphi)+C\eps(|\nabla u|^2,\varphi)
\]
for $\eps$ sufficiently small.

Using these estimates for $A_1,\,A_2,\,A_3$, we obtain
\begin{equation} \label{eq:lin}
((\Delta+1)u,  (\Delta+1)(u\varphi))
 \ge\tfrac34(|\Delta +1)u|^2,\varphi)-C\eps^2(u^2,\varphi)
-C\eps^2(|\nabla u|^2,\varphi).
\end{equation}

Next, using again integration by parts, $\sqrt{xy}\le x+y$ and~\eqref{eq:phi},
\begin{align*}
-(|\nabla u|^2,\varphi)
& =(\Delta u,u\varphi)+(u,\nabla u\cdot \nabla\varphi)
\\  &=((\Delta+1) u,u\varphi)-(u^2,\varphi)-\tfrac12(u^2,\Delta\varphi)
	\\  &\ge -{\textstyle\int}\{|(\Delta+1)u|\varphi^{1/2}\}\{|u|\varphi^{1/2}\}-(u^2,\varphi)-\tfrac12(u^2,\Delta\varphi)
\\  &\ge -(|(\Delta+1) u|^2,\varphi)-2(u^2,\varphi)-C\eps^2(u^2,\varphi).
\end{align*}
Substituting into~\eqref{eq:lin}, we obtain the desired result.
\end{proof}

\begin{cor} \label{cor:lin}
There exist constants $C>0$ and $\alpha>0$ such that
	\[
\frac{d}{dt}(u^2, \varphi) +(u^2,\varphi)
+
\alpha \big(\|u\|_{H^2_\varphi}^2+\|u\|_{L^4_\varphi}^4\big)
\le C(u^2,\varphi)-(u^4,\varphi)
\]
for all $\eps>0$ sufficiently small and $x_0\in\R^d$.
\end{cor}

\begin{proof}
By~\eqref{eq:L2} and Lemma~\ref{lem:lin},
\[
\frac{d}{dt}(u^2, \varphi) +(u^2,\varphi)+(u^4,\varphi)+(|(\Delta+1)u|^2,\varphi) \le C(u^2,\varphi)-(u^4,\varphi).
\]
By elliptic regularity for the Laplacian, 
for $\eps>0$ sufficiently small, there exists $\alpha>0$ independent of $x_0$  such that
\[
(u^4,\varphi)+(|(\Delta+1)u|^2,\varphi)\ge
	\alpha \big(\|u\|_{H^2_\varphi}^2+\|u\|_{L^4_\varphi}^4\big).
\]
The result follows.
\end{proof}

\begin{prop}[Dissipativity in $L^2_{\varphi}(\R^d)$] \label{prop:L2}
For $\eps>0$ sufficiently small, there is a constant $C_\eps>0$ such that
\[
	\|u(t)\|_{L^2_{\varphi}}^2+\int_{t_0}^t e^{-(t-s)} \Big(\|u(s)\|_{H^2_{\varphi}}^2  +\|u(s)\|_{L^4_{\varphi}}^4\Big)\,ds
\le
	C_\eps(e^{-(t-t_0)}\|u(t_0)\|_{L^2_{\varphi}}^2+1)
\]
for all  $x_0\in\R^d$, $0\le t_0\le t$.
\end{prop}

\begin{proof}
Change of variables gives $\int\varphi_{\eps,x_0}=\eps^{-d}\int\varphi_{1,0}$.
By Cauchy-Schwarz,
$(u^2,\varphi_{\eps,x_0})\le \|u^2\varphi_{\eps,x_0}^{1/2}\|_2
\|\varphi_{\eps,x_0}^{1/2}\|_2$, so squaring yields
\begin{equation} \label{eq:u4}
(u^4,\varphi_{\eps,x_0})\ge (u^2,\varphi_{\eps,x_0})^2 \|\varphi_{\eps,x_0}\|_1^{-1}
= C'\eps^d (u^2,\varphi_{\eps,x_0})^2 .
\end{equation}

Using this estimate, we obtain from Corollary~\ref{cor:lin} that
\begin{align*}
\frac d{dt}(u^2,\varphi)+(u^2,\varphi)
 + &  \alpha \Big(\|u\|_{H^2_{\varphi}}^2  +\|u\|_{L^4_{\varphi}}^4\Big)
  \le
C(u^2,\varphi)-C'\eps^d(u^2,\varphi)^2
\end{align*}
	Setting $v(t)=(u^2(t),\varphi)$ and integrating from $t_0$ to $t$, we arrive at
$$
	e^tv(t)
	-e^{t_0}v(t_0)
 +   \alpha \int_{t_0}^t e^s \Big(\|u(s)\|_{H^2_{\varphi}}^2  +\|u(s)\|_{L^4_{\varphi}}^4\Big)\,ds
\le
	C''\int_{t_0}^t e^s(Cv(s)-C'\eps^dv(s)^2)\,ds
$$
and the result follows.
\end{proof}

\begin{prop}[Dissipativity in $H^1_{\varphi}(\R^d)$] \label{prop:H1}
For $\eps>0$ sufficiently small, there is a constant $C_\eps>0$ such that
\[
	{\|u(t)\|}^2_{H^1_\varphi}
	+   \int_{t_0}^t e^{-(t-s)} {\|u(s)\|}_{H^3_\varphi}^2  \,ds
	\le C_\eps(e^{-(t-t_0)}\|u(t_0)\|^2_{H^1_\varphi}+1)
\]
for all $x_0\in\R^d$, $0\le t_0\le t$.
\end{prop}

\begin{proof}
We multiply equation~\eqref{eq:SH} by $-\divv(\varphi\nabla u)$ and integrate over $\R^d$.
Then
\begin{align*}
\frac12  \frac{d}{dt} \int|\nabla u|^2\varphi
	& =\int (\nabla u\cdot\nabla u_t)\varphi
=-\int \divv(\varphi\nabla u)u_t
\\ & =
\int \divv(\varphi\nabla u)(\Delta+1)^2u
-\lambda\int \divv(\varphi\nabla u)u
+\int \divv(\varphi\nabla u)u^3
\\ & = -\int ((\Delta+1)\nabla u) \cdot (\Delta+1)(\varphi\nabla u)
+\lambda\int |\nabla u|^2\varphi
-3\int u^2|\nabla u|^2\varphi.
\end{align*}
The linear terms are identical to those in~\eqref{eq:L2} with $u$ replaced by
$\nabla u$, so proceeding as before we obtain
	\[
		\frac{d}{dt} (|\nabla u|^2,\varphi)+
	(|\nabla u|^2,\varphi)+(u^2|\nabla u|^2,\varphi)+
	(|(\Delta+1)\nabla u|^2,\varphi)\le C(|\nabla u|^2,\varphi).
\]
By elliptic regularity,
	\[
		(u^2|\nabla u|^2,\varphi)+(|(\Delta+1)\nabla u|^2,\varphi)+(u,\varphi)\ge
	\alpha\|u\|_{H^3_\varphi}^2
\]
and hence
$$
\frac{d}{dt}(|\nabla u|^2, \varphi)
+(|\nabla u|^2, \varphi)
+\alpha\|u\|_{H_\varphi^3}^2
\le
 (u^2,\varphi)+C(|\nabla u|^2,\varphi).
$$

Integrating this inequality from $t_0$ to $t$, we obtain
\[
e^t(|\nabla u|^2(t),\varphi)
-e^{t_0}(|\nabla u|^2(t_0),\varphi)
 +   \int_{t_0}^t e^s \|u(s)\|_{H^3_{\varphi}}^2  \,ds
\le
C\int_{t_0}^t e^s \|u(s)\|_{H^2_{\varphi}}^2  \,ds.
\]
But
\[
	\int_{t_0}^t e^{-(t-s)} \|u(s)\|_{H^2_\varphi}^2  \,ds
\le
C_\eps(e^{-(t-t_0)}\|u(t_0)\|_{L^2_\varphi}^2+1)
\]
by Proposition~\ref{prop:L2}.
Hence
\[
	\|\nabla u(t)\|^2_{L^2_\varphi}
 +   \int_{t_0}^t e^{-(t-s)} \|u(s)\|_{H^3_\varphi}^2  \,ds
	\le C_\eps(e^{-(t-t_0)}\|u(t_0)\|^2_{H^1_\varphi}+1).
\]
The result follows from this combined with Proposition~\ref{prop:L2}.
\end{proof}

\begin{prop}[Dissipativity in $H^2_{\varphi}(\R^d)\cap L^4_{\varphi}(\R^d)$]
\label{prop:H2}
For $\eps>0$ sufficiently small, there is a constant $C_\eps>0$ such that
	\[
		\|u(t)\|^2_{H^2_{\varphi}}+\|u(t)\|^4_{L^4_{\varphi}}
	+\int_{t_0}^{t}e^{-(t-s)}\|u_t(s)\|^2_{L^2_{\varphi}}\,ds\le C_\eps(e^{-(t-t_0)}(\|u(t_0)\|^2_{H^2_{\varphi}}+\|u(t_0)\|^4_{L^4_{\varphi}})+1)
\]
for all $x_0\in\R^d$, $0\le t_0\le t$.
\end{prop}

\begin{proof}
Let $G(u)=\frac12 u^4-\lambda u^2$.
We multiply equation~\eqref{eq:SH} by $u_t\varphi$ and integrate over $\R^d$ to obtain
\[
(u_t^2,\varphi)+(u_t(\Delta+1)^2u,\varphi)+\tfrac12\tfrac{d}{dt}(G(u),\varphi)=0.
\]
Now, using integration by parts,
\[
\tfrac12\tfrac{d}{dt}(|(\Delta+1) u|^2,\varphi)=
(u_t(\Delta+1)^2u,\varphi)+
(u_t(\Delta+1)u,\Delta \varphi)+
2((\Delta+1)\nabla u\cdot\nabla\varphi,u_t).
\]
By~\eqref{eq:phi},
\[
|(u_t(\Delta+1)u,\Delta \varphi)|\le
(u_t^2,\Delta\varphi)+(|(\Delta+1)u|^2,\Delta\varphi)
\le C\eps^2 (u_t^2,\varphi)+C\eps^2(|(\Delta+1)u|^2,\varphi),
\]
and similarly
$|((\Delta+1)\nabla u\cdot\nabla\varphi,u_t)|
\le C\eps (u_t^2,\varphi)+C\eps(|(\Delta+1)\nabla u|^2,\varphi)$.
Hence, we arrive at
$$
(u_t^2,\varphi)+\frac d{dt}\Big((|(\Delta+1) u|^2,\varphi)+(G(u),\varphi)\Big)\le C\|u\|^2_{H^3_\varphi}
$$
and it follows that
\begin{align*}
	\frac d{dt}\Big((|(\Delta+1) u|^2,\varphi)+(G(u),\varphi)\Big)+\Big((|(\Delta+1) u|^2  ,\varphi) +  (G(u), & \varphi)  \Big)
	+  (u_t^2,\varphi)
	\\ &
\le C(\|u\|^2_{H^3_\varphi}+\|u\|^4_{L^4_\varphi}).
\end{align*}

	Integrating this inequality 
	from $t_0$ to $t$,
	\begin{align*}
		& \|(\Delta+1)u(t)\|^2_{L^2_{\varphi}}+  \|u(t)\|^4_{L^4_{\varphi}}
		 +\int_{t_0}^{t} e^{-(t-s)}\|u_t(s)\|^2_{L^2_{\varphi}}\,ds
		\\ & \le C  e^{-(t-t_0)}\Big(	\|(\Delta+1)u(t_0)\|^2_{L^2_{\varphi}}+ \|u(t_0)\|^4_{L^4_{\varphi}}\Big)
		+ C\int_{t_0}^t e^{-(t-s)}\Big(\|u(s)\|_{H_{\varphi}^3}^2+\|u(s)\|_{L_{\varphi}^4}^4\Big)\,ds.
		\end{align*}
Finally, we use Propositions~\ref{prop:L2} and~\ref{prop:H1} to control the integral on the right-hand side.
\end{proof}

\begin{prop}[Dissipativity of $u_t$ in $L^2_\varphi(\R^d)$] \label{prop:ut}
For $\eps>0$ sufficiently small, there is a constant $C_\eps>0$ such that
\[
 \|u_t(t)\|^2_{L^2_\varphi}\le C_\eps(e^{-(t-t_0)}\|u_t(t_0)\|^2_{L^2_\varphi}+1)
\]
for all $x_0\in\R^d$, $0\le t_0\le t$.
\end{prop}

\begin{proof}
Let $v=u_t$. Then
$
v_t+(\Delta+1)^2v = \lambda v-3u^2v.
$
Multiplying this equation by $v\varphi$ and using Lemma~\ref{lem:lin} (with $u$ replaced by $v$),
we obtain
\[
\tfrac{d}{dt}(v^2,\varphi)+(v^2,\varphi)\le C(v^2,\varphi).
\]
Integrating from $t_0$ to $t$, we obtain
\begin{align*}
 \|u_t(t)\|^2_{L^2_{\varphi}}
 \le   e^{-(t-t_0)}\|u_t(t_0)\|^2_{L^2_{\varphi}}
 + C\int_{t_0}^t e^{-(t-s)}\|u_t(s)\|_{L_{\varphi}^2}^2\,ds.
\end{align*}
Hence  by Proposition~\ref{prop:H2},
	\begin{equation} \label{eq:ut}
		\|u_t(t)\|^2_{L^2_{\varphi}}
	\le C_\eps(e^{-(t-t_0)}(\|u_t(t_0)\|^2_{L^2_\varphi}+\|u(t_0)\|^2_{H^2_{\varphi}}+\|u(t_0)\|^4_{L^4_{\varphi}})+1).
	\end{equation}

By Corollary~\ref{cor:lin},
\[
	\alpha(\|u\|^2_{H^2_\varphi} +\|u\|^4_{L^4_\varphi})
\le C(u^2,\varphi)-2(u_t\varphi^{1/2},u\varphi^{1/2})-(u^4,\varphi)
	\le C'\big\{(u^2,\varphi)+(u_t^2,\varphi)\big\}-(u^4,\varphi).
\]
By~\eqref{eq:u4} and completing the square, 
\(
\|u\|^2_{H^2_\varphi} +\|u\|^4_{L^4_\varphi}
	\le C''\eps^{-d}((u_t^2,\varphi)+1).
\)
Substituting this into~\eqref{eq:ut} yields the result.
\end{proof}

\subsection{Estimates in uniformly local spaces}

Define the uniformly local spaces $L^p_b(\R^d)$ with norm
$\|u\|_{L^p_b}=\sup_{x_0\in\R^d}\|u\mathbf{1}_{B_1(x_0)}\|_{L^p}$
and $H^k_b(\R^d)$ with norm
$\|u\|_{H^k_b}=\sup_{x_0\in\R^d}\|u\mathbf{1}_{B_1(x_0)}\|_{H^k}$.
By~\cite{Zelik03} (see also~\cite[Proposition~5.3]{MiranvilleZelik}),
for any $\eps>0$ there exists $C_\eps>1$ such that
\begin{equation} \label{eq:L2b}
	C_\eps^{-1}\sup_{x_0\in\R^d} (u^2,\varphi_{\eps,x_0}) \le
\|u\|_{L^2_b}^2\le C_\eps \sup_{x_0\in\R^d} (u^2,\varphi_{\eps,x_0}).
\end{equation}

\begin{prop}[Smoothing]
\label{prop:smooth}
	For all $t_1\in(0,1)$, there exists $C>0$ such that
\begin{itemize}
\item[(a)] $\|u(t)\|^2_{H^2_b}+\|u(t)\|^4_{L^4_b} \le C(e^{-t}\|u_0\|^2_{L^2_b}+1)$ and
\item[(b)] $\|u_t(t)\|^2_{L^2_b} \le C(e^{-t}\|u_0\|^2_{L^2_b}+1)$
\end{itemize}
	for all $t\ge t_1$.
	\end{prop}

	\begin{proof}
		Fix $\eps>0$ sufficiently small
		and let $\varphi=\varphi_{\eps,x_0}$.

\vspace{1ex} \noindent(a)
		By Proposition~\ref{prop:L2} and estimate~\eqref{eq:L2b},
		there is a constant $C>0$ such that
\[
	\int_{0}^{t_1} \Big(\|u(s)\|_{H^2_{\varphi}}^2  +\|u(s)\|_{L^4_{\varphi}}^4\Big)\,ds
\le
C (\|u_0\|_{L^2_b}^2+1)
\]
		for all $x_0\in\R^d$.
Hence there exists $t_0\in(0,t_1)$ depending on $u$ and $x_0$ such that
$$
		\|u(t_0)\|_{H^2_{\varphi}}^2+\|u(t_0)\|^4_{L^4_{\varphi}}\le
		Ct_1^{-1}(\|u_0\|^2_{L^2_b}+1).
$$
By Proposition~\ref{prop:H2},
$$
		\|u(t)\|^2_{H^2_{\varphi}}+\|u(t)\|^4_{L^4_\varphi}\le C(e^{-t}(\|u(t_0)\|_{H^2_{\varphi}}^2+\|u(t_0)\|^4_{L^4_{\varphi}}) +1),
\quad t\ge t_0.
$$
Combining these estimates, we obtain that
$$
		\|u(t)\|^2_{H^2_{\varphi}}
+\|u(t)\|^4_{L^4_\varphi}
		\le Ct_1^{-1}(e^{-t}
		\|u_0\|^2_{L^2_b}+1)
$$
for all $t\ge t_0$ and hence for all $t\ge t_1$.
		Now apply~\eqref{eq:L2b} once more.

\vspace{1ex} \noindent (b)
By Proposition~\ref{prop:H2}, part~(a) and estimate~\eqref{eq:L2b},
\[
\int_{\frac12 t_1}^{t_1}\|u_t(s)\|_{L^2_\varphi}^2\le 
		C(\|u(\tfrac12 t_1)\|_{H^2_\varphi}^2+\|u(\tfrac12 t_1)\|_{L^4_\varphi}^4+1)
\le C(\|u_0\|_{L^2_b}^2+1).
\]
Hence there exists $t_0\in(0,t_1)$ depending on $u$ and $x_0$ such that
\[
\|u_t(t_0)\|_{L^2_\varphi}^2\le Ct_1^{-1}(\|u_0\|_{L^2_b}^2+1).
\]
Now combine this with Proposition~\ref{prop:ut} and apply estimate~\eqref{eq:L2b}.
	\end{proof}

\subsection{Estimates in $L^\infty(\R^d)$}
\label{sec:8}

Classical Sobolev embedding theorems carry over to uniformly local spaces.
For example, suppose that $d<2k$. Then $H^k(\R^d)$ is embedded in $L^\infty(\R^d)$. Hence there is a constant $C>0$ independent of $x_0$ such that
$\|u\mathbf{1}_{B_1(x_0)}\|_{\infty}\le C\|u\mathbf{1}_{B_1(x_0)}\|_{H^k}$.
Taking the supremum over $x_0$, we obtain
$\|u\|_{\infty}\le C\|u\|_{H^k_b}$.

Set $F(u)=-\Delta^2u+\lambda u-u^3$.

\begin{lemma} \label{lem:embed}
For $1\le d\le 9$, there is a constant $C>0$ such that
\[
\|u\|_{H_b^4}\le C(\|u\|_{H_b^2}+\|F(u)\|_{L_b^2})
	\quad\text{for all $u\in H_b^2$ such that $F(u)\in L^2_b$.}
\]
\end{lemma}

\begin{proof}
First we show (without any restriction on $d$) that
\begin{equation} \label{eq:emb}
\|u\|_{H_b^3}\le C(\|u\|_{H_b^2}+\|F(u)\|_{L_b^2}).
\end{equation}

Let $g(u)=\lambda u-u^3$.
	Multiplying the equation $F(u)=-\Delta^2u+g(u)$ by $\divv(\varphi\nabla u)$ and integrating over $\R^d$, we obtain
\begin{align*}
\int |\nabla\Delta u|^2\varphi & = \int \nabla\Delta u\cdot \nabla\divv(\varphi\nabla u)-
\int \nabla\Delta u\cdot \big\{ \nabla\varphi\Delta u
+\nabla(\nabla\varphi\cdot\nabla u) \big\}
\\ & = -\int \Delta^2 u \divv(\varphi\nabla u)-\int \nabla\Delta u\cdot \big\{ \nabla\varphi\Delta u
+\nabla(\nabla\varphi\cdot\nabla u) \big\}
\\ & =\int \varphi g'(u)|\nabla u|^2 
+\int F(u)(\varphi \Delta u+\nabla\varphi\cdot\nabla u)
\\ & \qquad \qquad \qquad\qquad -\int \nabla\Delta u\cdot \big\{ \nabla\varphi\Delta u
+\nabla(\nabla\varphi\cdot\nabla u) \big\}.
\end{align*}
Using the estimates $g'(u)=\lambda-3u^2+\lambda\le\lambda_0$ and~\eqref{eq:phi},
\begin{align*}
\|\nabla\Delta u\|^2_{L^2_\varphi}
 & \le \lambda_0\int |\nabla u|^2\varphi  + \int F(u)^2\varphi  + \int |\Delta u|^2\varphi + \eps\|u\|^2_{H^3_\varphi}
\\ & \le C\|u\|^2_{H^2_\varphi}+ \|F(u)\|^2_{L^2_\varphi}+\eps\|u\|^2_{H^3_\varphi}.
\end{align*}
Hence we obtain~\eqref{eq:emb} by
choosing $\eps$ sufficiently small and applying~\eqref{eq:L2b}.

Now suppose that $u\in H^2_b$ and
$F(u)\in L^2_b$.
By~\eqref{eq:emb} and the Sobolev embedding theorem (since $d\le 9$), $u\in H^3_b\subset L^6_b$. Hence $u^3\in L^2_b$ and it follows that $g(u)\in L^2_b$.
But this means that $\Delta^2u=g(u)-F(u)\in L^2_b$
and hence $u\in H^4_b$.
\end{proof}

\begin{pfof}{Theorem~\ref{thm:app} for $d\le 7$.}
The inequality $\|u\|_{L_b^2}\le C\|u\|_{\infty}$ holds for all $d\ge1$.
By the Sobolev embedding theorem,
$\|u\|_{\infty} \le C\|u\|_{H^4_b}$ for $d\le 7$. Hence
by Proposition~\ref{prop:smooth} and Lemma~\ref{lem:embed}, for all $t_1\in(0,1)$, there exists $C>0$ such that
\begin{align*}
\|u(t)\|_{\infty} & \le C\|u(t)\|_{H^4_b}
\le C(\|u(t)\|_{H^2_b}+\|F(u(t))\|_{L^2_b})
= C(\|u(t)\|_{H^2_b}+\|u_t(t)\|_{L^2_b})
\\ & \le C(e^{-t}\|u_0\|_{L^2_b}+1)
\le C(e^{-t}\|u_0\|_{\infty}+1)
\end{align*}
for all $t\ge t_1$.
Combining this with local well-posedness in $L^\infty(\R^d)$, there exists 
$C>0$ such that
$$
\|u(t)\|_{\infty} 
\le C(e^{-t}\|u_0\|_{\infty}+1)
\quad\text{for all $t\ge0$}
$$
as required.
\end{pfof}

The result for $d\le 9$ can be obtained by bootstrapping.
The argument above shows that if $u_0\in L^\infty$
then $u(t)\in H^4_b$ which is embedded in $L^{18}_b$ so $g(u(t))\in L^6_b$.
Now we treat equation~\eqref{eq:SH} as a linear equation with right-hand side $g(u)\in L^\infty([0,\infty), L^6_b(\R^d))$. Parabolic regularity then implies that $u\in L^\infty([0,\infty), W^{k,6}_b(\R^d))$ for any $k<4$ and this is embedded in $L^\infty([0,\infty)\times \R^d))$.

\end{document}